\documentclass{article}

\usepackage{amsmath,amsfonts,amssymb,a4,enumerate,epsfig,theorem,psfrag}
\usepackage[latin1]{inputenc}

\newcommand{\inv}{\textrm{\rm inv}}

\newcommand{\tub}{\textrm{\rm tub}}

\newcommand{\bbb}{{b}}

\newcommand{\dom}{{\operatorname{Dom}}}
\newcommand{\dist}{{\operatorname{dist}}}
\newcommand{\sign}{{\operatorname{sign}}}
\newcommand{\trace}{{\operatorname{tr}}}
\newcommand{\transpose}{{\ast}}

\newcommand{\lolo}{{L}}

\newcommand{\NN}{{\mathbb{N}}}
\newcommand{\BB}{{\mathbb{B}}}
\newcommand{\ZZ}{{\mathbb{Z}}}

\newcommand{\RR}{{\mathbb{R}}}

\newcommand{\Mu}{{M}}
\newcommand{\Aa}{{\cal{A}}}
\newcommand{\Bb}{{\cal{B}}}
\newcommand{\Cc}{{\cal{C}}}
\newcommand{\Dd}{{\cal{D}}}
\newcommand{\Ee}{{\cal{E}}}
\newcommand{\cE}{{\cal{E}}}

\newcommand{\Ii}{{\cal{I}}}
\newcommand{\Jj}{{\cal{J}}}
\newcommand{\Kk}{{\cal{K}}}
\newcommand{\cK}{{\cal{K}}}

\newcommand{\Oo}{{\cal{O}}}

\newcommand{\cP}{{\cal{P}}}
\newcommand{\cQ}{{\cal{Q}}}

\newcommand{\calS}{{\cal{S}}}
\newcommand{\Tt}{{\cal{T}}}
\newcommand{\cT}{{\cal{T}}}
\newcommand{\Uu}{{\cal{U}}}

\newcommand{\cX}{{\cal{X}}}

\newcommand{\cN}{{\cal{N}}}
\newcommand{\cM}{{\cal{M}}}
\newcommand{\Xx}{{\cal{X}}}

\newcommand{\ILa}{\Tt_\Lambda}
\newcommand{\OLa}{\Oo_\Lambda}
\newcommand{\OSa}{\Oo_\Sigma}

\newcommand{\ILai}{\Tt_{\Lambda_i}}

\newcommand{\JLa}{\bar\Jj_\Lambda}

\newcommand{\JLao}{\Jj_\Lambda}

\newcommand{\proof}{{\noindent\bf Proof: }}

\def\qed{\unskip\nobreak\hfil\penalty50\hskip1.75em\null\nobreak\hfil
$\blacksquare$ {\parfillskip=0pt \finalhyphendemerits=0 \par}\medbreak}

\newcommand\capsize{\relax}
\newcommand\nobf{\noindent\bf}

\newcommand\grad{\operatorname{grad}}
\newcommand\diag{\operatorname{diag}}

\newcommand\interior{\operatorname{int}}

\newtheorem{theo}{Theorem}
\newtheorem{lemma}{Lemma}[section]
\newtheorem{prop}[lemma]{Proposition}
\newtheorem{coro}[lemma]{Corollary}

\newtheorem{fact}[lemma]{Fact}

\title{Dynamics of the symmetric eigenvalue problem with shift strategies}
\author{Ricardo S. Leite, Nicolau C. Saldanha and Carlos Tomei}

\begin{document}

\maketitle

\begin{abstract}
A common algorithm for the computation of eigenvalues
of real symmetric tridiagonal matrices 
is the iteration of certain special maps $F_\sigma$
called \textit{shifted $QR$ steps}.
Such maps preserve spectrum and
a natural common domain is $\ILa$,
the manifold of  real symmetric tridiagonal matrices
conjugate to the diagonal matrix $\Lambda$.
More precisely, a (generic) \textit{shift} $s \in \RR$
defines a map $F_s: \ILa \to \ILa$.
A strategy $\sigma: \ILa \to \RR$ specifies the shift to be applied at $T$
so that $F_\sigma(T) = F_{\sigma(T)}(T)$.
Good shift strategies should lead to fast \textit{deflation}:
some off-diagonal coordinate tends to zero,
allowing for reducing of the problem to submatrices.
For topological reasons,
continuous shift strategies do not obtain fast deflation;
many standard strategies are indeed discontinuous.
Practical implementation only gives rise systematically
to \textit{bottom deflation},
convergence to zero of the lowest off-diagonal entry $\bbb(T)$.
For most shift strategies, convergence to zero of $\bbb(T)$ is
cubic, $|\bbb(F_\sigma(T))| = \Theta(|\bbb(T)|^k)$ for $k = 3$.
The existence of arithmetic progressions in the spectrum of $T$
sometimes implies instead quadratic convergence, $k = 2$.
The complete integrability of the Toda lattice
and the dynamics at non-smooth points
are central to our discussion.
The text does not assume knowledge of numerical linear algebra.
\end{abstract}


\medbreak


\smallbreak

{\noindent\bf MSC2010-class:} 37N30, 37C05, 65F15.

\section{Introduction}

In this paper, we study some subtle dynamical aspects
of a class of numerical algorithms
for eigenvalues of real symmetric matrices.
This includes the classic
inverse iteration with different shift strategies,
among which Rayleigh and Wilkinson shifts.
We do not assume previous knowledge of numerical linear algebra.

Numerical analysts are familiar with \textit{tridiagonalization},
the fact that given a real symmetric matrix $S$
it is easy to obtain another isospectral
matrix $T$ which is tridiagonal: $(T)_{ij} = 0$ whenever $|i-j| > 1$.
For matrices of order approximately between $20$ and $1000$,
it pays to first tridiagonalize and then work in the vector space $\cT$
of real symmetric tridiagonal matrices.
Let $\Lambda = \diag(\lambda_1 < \lambda_2 < \cdots < \lambda_n)$
be a diagonal matrix with simple spectrum:
it turns out that the set $\ILa \subset \cT$
of tridiagonal matrices isospectral with $\Lambda$
is a connected compact smooth oriented manifold (\cite{Tomei}, \cite{LST1}).
The algorithms under consideration are defined by iteration
of some easily computable map $F: \ILa \to \ILa$:
given $T \in \ILa$ we consider the sequence $(F^k(T_0))$.
For relevant maps $F$,
diagonal matrices in $\ILa$ are fixed points of $F$.

Let $\cE \subset O(n)$ be the group of
real orthogonal diagonal matrices,
so that for $E \in \cE$, $(E)_{ii} = \pm 1$.
For each $E \in \cE$, the map $\eta: \ILa \to \ILa$, $\eta(T) = ETE$,
is an involutive diffeomorphism of $\ILa$:
its effect on $T \in \ILa$ is to change signs
of some subdiagonal entries $(T)_{i+1,i}$.
Numerical analysts, familiar with this simple fact,
often drop signs of subdiagonal entries.
We shall not do likewise for we are often interested
in smoothness issues.
Again, relevant maps will be ($\cE$-)\textit{equivariant},
in the sense that $F \circ \eta = \eta \circ F$ for all $\eta$.

\begin{figure}[ht]
\centerline{\epsfig{height=65mm,file=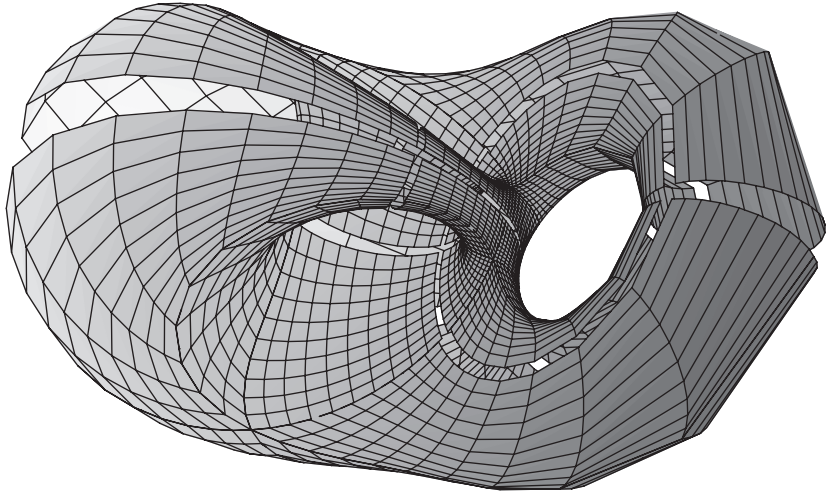}}
\caption{\capsize The cell decomposition
of $\ILa$ for $\Lambda = \diag(4,5,7)$}
\label{fig:bitorus}
\end{figure}

Consistently with the involutions above,
signs of subdiagonal entries
induce a cell decomposition of $\ILa$.
The $0$-cells are the $n!$ diagonal matrices
and the top dimensional $(n-1)$-cells
turn out to be $2^{n-1}$ permutohedra
(polytopes equivalent to the convex hull of
the $n!$ points of $\RR^n$ obtained
by permuting $n$ fixed distinct real numbers).
For $n = 3$, the manifold $\ILa$ is a bitorus
which can be obtained by gluing four hexagons
along six circles (see Figure \ref{fig:bitorus}).

We consider that iteration of the map $F$
has accomplished its job when one subdiagonal entry
$(F^k(T))_{i+1,i}$ has absolute value smaller than some prescribed tolerance.
In terms of the cell decomposition, 
we are done when we hit (a thin neighborhood of) a lower dimensional cell,
or, in the numerical jargon,
the sequence $(F^k(T))$ undergoes \textit{deflation}.
Notice that if $(T)_{i+1,i} = 0$
then the matrix $T$ \textit{splits} as $T = T^a \oplus T^b$
where the tridiagonal submatrices $T^a$ and $T^b$
have orders $i$ and $n-i$, respectively.
Pragmatically, if $(T)_{i+1,i} \approx 0$ then
the spectrum of $T$ is approximately the disjoint union
of the spectra of $T^a$ and $T^b$, which are easier to compute.

Ideally, deflation should happen approximately in the middle
so that each subproblem has order
approximately half of the original one.
Unfortunately, it is not known how to implement
easily computable iterations with this property.
Usually the sequence $(F^k(T))$ undergoes \textit{bottom deflation}:
\[ \lim_{k \to +\infty} \bbb(F^k(T)) = 0; \qquad
\bbb(T) = (T)_{n,n-1}. \]
Geometrically, we approach one of the $n$ \textit{deflation sets}
$\Dd^i_{\Lambda,0} \subset \ILa$ defined by
$(T)_{n,n} = \lambda_i$, $\bbb(T) = 0$.

Notice that removing the $n$-th row and column
obtains a diffeomorphism:
\[ \Dd^i_{\Lambda,0} \approx \cT_{\Lambda_i}, \qquad
\Lambda_i = \diag(\lambda_1, \ldots, \lambda_{i-1},
\lambda_{i+1}, \ldots, \lambda_n); \]
in particular, $\Dd^i_{\Lambda,0}$ is connected.
In Figure \ref{fig:bitorus}, the submanifolds $\Dd^i_{\Lambda,0}$
are three of the six (removed) circles.
It turns out (Proposition \ref{prop:deflation}) that,
for sufficiently small $\epsilon > 0$,
the closed set $\Dd_{\Lambda,\epsilon} \subset \ILa$
defined by $|\bbb(T)| \le \epsilon$ has $n$ connected components
$\Dd^i_{\Lambda,\epsilon}$ which are 
closed tubular neighborhoods of $\Dd^i_{\Lambda,0}$.

\smallbreak

As algorithms for eigenvalue computation,
continuous maps $F$ are problematic.

\begin{theo}
\label{theo:connect}
Let $F: \ILa \to \ILa$ be a continuous $\cE$-equivariant map
such that every diagonal matrix in $\ILa$ is a fixed point of $F$.
\begin{enumerate}[(a)]
\item{The map $F$ is surjective.}
\item{If there exist disjoint compact sets $\cK_i \supset \Dd^i_{\Lambda,0}$
with $F(\cK_i) \subset \interior(\cK_i)$
then there exists $T \in \ILa$ for which
the sequence $(F^k(T))$ does not undergo bottom deflation.}
\end{enumerate}
\end{theo}

Item (a) already makes $F$ unpromising as an algorithm:
for any $k$, the iterate $F^k$ is surjective
and, given $k$, there exists $T$ such that $F^k(T)$ is far from deflation.
The additional hypothesis in item (b),
which, as we shall see, holds for many algorithms,
makes $F$ even less desirable.
The proof of this result uses methods very different
from the rest of the paper and is left for the Appendix.

These phenomena lead numerical analysts to consider discontinuous maps $F$.
Among the standard algorithms to compute eigenvalues
of matrices in $\Tt$ are $QR$ steps
with different shift strategies: Rayleigh and Wilkinson are familiar examples
(excellent references are \cite{Wilkinson}, \cite{Demmel}, \cite{Parlett}).
Recall that Rayleigh's strategy $\rho$ is continuous and is known
to have the unfortunate property (b) that there exists
a matrix $T$ for which $(F_\rho^k(T))$ does not undergo bottom deflation;
Wilkinson, on the other hand, is discontinuous.
In this paper, we consider a more general context: we define simple shift
strategies, which include the examples above and more.


More precisely, given a matrix $T \in \Tt$ and $s \in \RR$,
write $T - sI = Q R$, if possible, for an orthogonal matrix $Q$
and an upper triangular matrix $R$ with positive diagonal entries.
A \textit{shifted $QR$ step} is $\Phi(T,s) = Q^\transpose T Q$.
As is well known, shifted $QR$ steps preserve spectrum and shape.
A function $\sigma: \ILa \to \RR$ is ($\cE$-)\textit{invariant}
if $\sigma(\eta(T)) = \sigma(ETE) = \sigma(T)$
for all $T \in \ILa$ and all $E \in \cE$.
A \textit{simple shift strategy} is
an invariant function $\sigma: \ILa \to \RR$
satisfying the following condition:
there exists $C_\sigma > 0$ such that for all $T \in \ILa$
there is an eigenvalue $\lambda_i$ with $| \sigma(T) - \lambda_i | \le
C_\sigma |\bbb(T)|$.

For technical reasons, we prefer the \textit{signed} step   
$F_s(T) = \Phi_\star(T,s) = Q_\star^\transpose T Q_\star$,
where now $T - sI = Q_\star R_\star$,
the orthogonal matrix $Q_\star$ has positive determinant and
only the first $n-1$ diagonal entries of the upper triangular matrix $R_\star$
are required to be positive.
As we shall see, the signed step is smoothly defined on a larger domain,
and convergence issues for both kinds of step iterations
are essentially equivalent.


Simple shift strategies prescribe shifts:
set $F_\sigma(T)= F_{\sigma(T)}(T)$.
It turns out that $F_\sigma$ is a well-defined
(but usually discontinuous) equivariant map from
$\ILa$ (or some very large subset thereof) to $\ILa$.

An important question in practice is estimating
the rate of deflation, i.e.,
the rate of convergence to zero of the sequence $\bbb(F_\sigma^k(T))$.
Numerical evidence indicates that deflation is often cubic,
in the sense that there is a constant $C$ such that
$|\bbb(F_\sigma^{k+1}(T))| \le C |\bbb(F_\sigma^k(T))|^3$ for large $k$.

Consider the \textit{singular support}
$\calS_\sigma \subset \ILa$ of a shift strategy $\sigma$,
the minimal closed subset of $\ILa$
on whose complement $\sigma$ is smooth.
Away from the singular support $\calS_\sigma$, squeezing is cubic.

\begin{theo}
\label{theo:squeeze}
For $\epsilon > 0$ small enough,
each deflation neighborhood $\Dd^i_{\Lambda,\epsilon}$
is invariant under $F_\sigma$.
There exists $C > 0$ such that,
for all $T \in \Dd_{\Lambda,\epsilon}$,
$|\bbb(F_\sigma(T))| \le C |\bbb(T)|^2$.
Also, given a compact set $\Kk\subset \Dd_{\Lambda,\epsilon}$
disjoint from $\calS_\sigma \cap \Dd_{\Lambda,0}$,
there exists $C_\Kk > 0$ such that,
for all $T \in \Kk$,
$|\bbb(F_\sigma(T))| \le C_{\Kk} |\bbb(T)|^3$.
\end{theo}

Although the tubular neighborhoods $\Dd^i_{\Lambda,\epsilon}$
are invariant under $F_\sigma$,
it is not true in general
that $F_\sigma^k(T)$ belongs to  $\Dd^i_{\Lambda,\epsilon}$
for sufficiently large $k$:
this is true, however, for the important example
of Wilkinson's shift.

\begin{figure}[ht]
\psfrag{X0}{$\Xx_0$}
\psfrag{X1p}{$\Xx_1^{(+)}$}
\psfrag{X1m}{$\Xx_1^{(-)}$}
\psfrag{X2pp}{$\Xx_2^{(++)}$}
\psfrag{X2pm}{$\Xx_2^{(+-)}$}
\psfrag{X2mp}{$\Xx_2^{(-+)}$}
\psfrag{X2mm}{$\Xx_2^{(--)}$}
\centerline{\epsfig{height=25mm,file=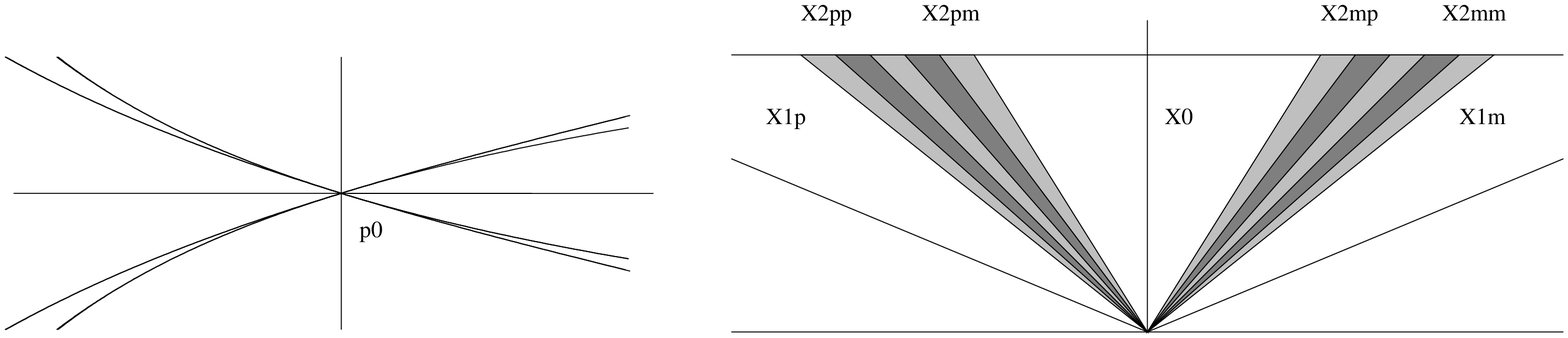}}
\caption{\capsize The set $\cX$}
\label{fig:xaxbxc}
\end{figure}

For Rayleigh's shift, it is well known that convergence
(when it happens) is always cubic;
this is a corollary of Theorem \ref{theo:squeeze}.
Cubic convergence does not hold in general for Wilkinson's strategy.
In \cite{LST2}, for $\Lambda= \diag(-1,0,1)$,
we construct a Cantor-like set $\cX \subset \ILa$
of unreduced initial conditions
for which the rate of convergence is strictly quadratic.
Sequences starting at $\cX$ converge to a reduced
matrix which is not diagonal.
A part of the set $\cX$ is shown (true to scale)
on the left part of Figure \ref{fig:xaxbxc};
the Cantor-like aspect is invisible:
the cross section of each of the four visible ``curves''
really consists of a tiny Cantor set,
far smaller than the resolution of the picture.
This is consistent with the fact that
the Hausdorff dimension of $\cX$ is $1$.
The set $\cX$ is the intersection of thinner and thinner wedges;
in the right half we show schematically
three generations of such wedges.
The central vertical line is a step-like discontinuity;
the inverse image of the largest wedge $\Xx_0$ is
$\Xx_1^{(+)} \cup \Xx_1^{(-)}$
and the inverse image of \textit{that} is
$\Xx_2^{(++)} \cup \Xx_2^{(+-)} \cup
\Xx_2^{(-+)} \cup \Xx_2^{(--)}$.

As numerical analysts know,
shift strategies usually define sequences of matrices which, asymptotically,
not only isolate an eigenvalue at the $(n,n)$ position but also isolate,
at a slower rate, a second eigenvalue at the $(n-1,n-1)$ position.
This does not happen for the example above
where $(F_\sigma^k(T))_{n,n}$ tends to the center
of a three-term arithmetic progression of eigenvalues and
$(F_\sigma^k(T))_{n-1,n-2}$ stays bounded away from zero.

A matrix $T \in \Tt$ with simple spectrum is {\it a.p.\ free} if it does not
have three eigenvalues in arithmetic progression and {\it a.p.}\ otherwise;
in particular, generic spectra are a.p.\ free.
In this case, the situation is very nice:
cubic convergence is essentially uniform on $\ILa$.
This condition is reminiscent of the Sternberg's resonance hypothesis
for normal forms (\cite{Sternberg}).


\begin{theo}
\label{theo:big}
Let $\Lambda$ be an a.p.\ free matrix and
$\sigma$ a shift strategy
for which diagonal matrices do not belong to $\calS_\sigma$.
Then there exist $\epsilon > 0$, $C >0$ and $K > 0$ such that:
\begin{enumerate}[(a)]
\item{the deflation neighborhood $\Dd_{\Lambda,\epsilon}$
is invariant under $F_\sigma$;}
\item{for any $T \in \Dd_{\Lambda,\epsilon}$,
the sequence $(F_\sigma^k(T))$ converges to a diagonal matrix
and the set of positive integers $k$ for which
$|\bbb(F_\sigma^{k+1}(T))| > C |\bbb(F_\sigma^{k}(T))|^3$
has at most $K$ elements.}
\end{enumerate}
\end{theo}

Still, the finite set of points in which the cubic estimate does not
hold may occur arbitrarily late along the sequence $(F_\sigma^{k}(T))$.

An a.p.\ matrix is {\it strong a.p.}\ if it contains three consecutive
eigenvalues in arithmetic progression and {\it weak a.p.}\ otherwise.
Under very mild additional hypothesis, $\bbb(T)$ converges to zero at a cubic
rate also for weak a.p.\ matrices.
Let $\Cc_{\Lambda,0} \subset \ILa$ be the set of matrices $T$
for which $(T)_{n,n-1} = (T)_{n-1,n-2} = 0$.

\begin{theo}
\label{theo:bigg}
Let $\Lambda$ be a weak a.p.\ matrix and
$\sigma: \ILa \to \RR$ a shift strategy
for which $\Cc_{\Lambda,0}$ and $\calS_\sigma$ are disjoint.
Then there exists $\epsilon > 0$ such that
the deflation neighborhood  $\Dd_{\Lambda,\epsilon}$
is invariant under $F_\sigma$
and, for all unreduced $T \in \Dd_{\Lambda,\epsilon}$,
the sequence $(\bbb(F_\sigma^k(T)))$ converges to zero
at a rate which is at least cubic.
More precisely,
for each unreduced $T \in \Dd_{\Lambda,\epsilon}$
there exist $C_T, K_T > 0$ such that,
for all $k > K_T$,
we have $|\bbb(F_\sigma^{k+1}(T))| \le C_T |\bbb(F_\sigma^k(T))|^3$.
\end{theo}

In particular, the convergence of Wilkinson's strategy
is cubic for weak a.p.\ matrices.
However, uniformity in the sense of Theorem \ref{theo:big}
is not guaranteed and
the constants $C_T$ and $K_T$ depend on $T$.
As in the case of the spectrum $\{-1,0,1\}$,
we conjecture that if $\Lambda$ is strong a.p.\ then
there exists $\cX \subset \ILa$
of Hausdorff codimension $1$ of initial conditions $T$ for which the rate of
convergence is strictly quadratic.

\smallskip


The celebrated integrability of the Toda lattice
(\cite{Flaschka}, \cite{Moser})
on unreduced tridiagonal matrices 
manifests itself in several ways along the paper:
it provided ample inspiration
but the paper strives to be self-contained.
For starters, the steps $F_s$, $s \in \RR$,
commute in their natural domains (Proposition \ref{prop:commute}).
Norming constants (as in \cite{Moser})
provide angle variables for which steps $F_s$
are translations.
Unfortunately, these angle variables break down
(as they must!) for reduced matrices $T \in \ILa$.
Since $(F^k_\sigma(T))$ approaches reduced matrices
we prefer to introduce other coordinate systems
which extend smoothly to such points.
\textit{Bidiagonal coordinates}, defined in \cite{LST1},
consist of very explicit charts on the manifold $\ILa$.
They are used in \cite{LST1}
to prove the cubic convergence of Rayleigh's shift
and in the unpublished manuscript \cite{LST3}
to prove some of the results presented here for Wilkinson's shift.
In Section 4, instead, we introduce \textit{tubular coordinates}
on the tubular neighborhoods $\Dd_{\Lambda,\epsilon}^i$:
steps $F_s$ within these sets are given by a very simple formula
(Corollary \ref{coro:tubularQR}).


The (signed) steps $F_\sigma$ are smooth whenever the shift strategy is, i.e.,
for $T \in \Dd_{\Lambda,\epsilon}^i \smallsetminus \calS_\sigma$
(``unsigned steps'' would not be smooth on limit points).
At matrices $T_0 \in \Dd^i_{\Lambda,0}$ on which $F_\sigma$ is smooth,
the map $T \mapsto \bbb(F_\sigma(T))$ has zero gradient.
The symmetry of the shift strategy
yields a cubic Taylor expansion and therefore an estimate
$|\bbb(F_\sigma(T))| \le C |\bbb(T)|^3$, settling Theorem \ref{theo:squeeze}.

\textit{Height functions} $H: \Dd_{\Lambda,\epsilon}^i \to \RR$
(similar to Lyapunov functions)
are used for further study of the sequence $(F_\sigma^k(T))$
in the a.p.\ free case.
More precisely, for steps $s$ near $\lambda_i$,
$H_i(F_s(T))> H_i(T)$ provided
$T \in \Dd^i_{\Lambda,\epsilon}$ is not diagonal:
this is another manifestation of the Toda dynamics.
Theorem \ref{theo:big} then follows by
a compactness argument bounding the number of iterations
for which $F_\sigma^k(T)$ stays close
to the singular support ${\calS}_\sigma$.

For a.p. spectra the situation is subtler,
as can be seen from the example in \cite{LST2}
and Figure \ref{fig:xaxbxc}.
On the other hand, Theorem \ref{theo:bigg} tells us that
the weak a.p.\ hypothesis
together with an appropriate smoothness condition
guarantee cubic convergence.

\smallskip

In Section 2 we list the basic properties of the signed shifted $QR$ step
on the manifold $\ILa$.
Simple shift strategies are introduced in Section 3,
and the standard examples are shown to satisfy the definition.
We define the deflation set $\Dd_{\Lambda,0}$
and neighborhood $\Dd_{\Lambda,\epsilon}$
in Section 4 and then set up tubular coordinates.
The local theory of steps $F_s$ near $\Dd_{\Lambda,0}$
and the proof Theorem \ref{theo:squeeze} are presented in Section 5.
In Section 6 we construct the height functions $H$
and then prove Theorem \ref{theo:big}.
The convergence properties for a.p.\ matrices
in Theorem \ref{theo:bigg} are proved in Section 7.
We present in Section 8 two counterexamples
to natural but incorrect strengthenings
of Theorems \ref{theo:big} and \ref{theo:bigg}.
Finally, the Appendix is dedicated to Theorem \ref{theo:connect}.


The authors are very grateful for the abundant contributions 
of several readers of this work and its previous versions.
The authors acknowledge support from CNPq, CAPES, INCT-Mat and FAPERJ.


\section{The manifold $\ILa$ and shifted steps $F_s$}

Let $\Tt$ denote the real vector space of $n \times n$ real, symmetric,
tridiagonal matrices endowed with the norm $\|T\|^2 = \trace(T^2)$.
For $T \in \Tt$, the \textit{subdiagonal entries}
of $T$ are $(T)_{i+1,i}$ for $i=1,\ldots,n-1$.
The lowest subdiagonal entry of $T$ is $\bbb(T)= (T)_{n,n-1}$.

As usual, let $SO(n)$ denote the set of orthogonal matrices
with determinant equal to $1$.
Let $\Lambda$ be a real diagonal matrix
with simple eigenvalues $\lambda_1 < \cdots < \lambda_n$.
Define the \textit{isospectral manifold}
\[ \ILa = \{ Q^\ast \Lambda Q, Q \in SO(n)\} \cap \Tt, \]
the set of matrices in $\Tt$ similar to $\Lambda$.
The set $\ILa \subset \Tt$ is a real smooth manifold
(\cite{Tomei}; \cite{LST1} describes an explicit atlas of $\ILa$).


For a matrix $M$, the $QR$ factorization is $M = QR$
for an orthogonal matrix $Q$ and
an upper triangular matrix $R$ with positive diagonal.
The {\it $Q_\star R_\star$ factorization}, instead,
is $M = Q_\star R_\star$, for $Q_\star \in SO(n)$ and
$R_\star$ an upper triangular matrix
with $(R_\star)_{i,i} > 0$, $i=1,\ldots,n-1$.
A real $n \times n$ matrix $M$ is {\it almost invertible}
if its first $n-1$ columns are linearly independent:
notice that almost invertible matrices are dense
within $n \times n$ matrices and form an open set.
The diagonal matrix $E_{n-1}$ is such that $(E_{n-1})_{i,i}$
is $1$ for $i < n$ and $-1$ for $i=n$.

\begin{prop}
\label{prop:existence}
An almost invertible real matrix $M$ admits
a unique $Q_\star R_\star$ factorization,
with $Q_\star$ and $R_\star$ depending smoothly
on $M$.
If $M$ is invertible, it admits unique (smooth)
factorizations $M = QR = Q_\star R_\star$.
If $\det M >0$, the factorizations are equal, i.e.,
$Q= Q_\star$ and  $R = R_\star$.
If $\det M <0$, $Q= Q_\star E_{n-1}$ and $R = E_{n-1} R_\star$.
If $\det M =0$, $(R_\star)_{n,n}=0$.
\end{prop}


\proof
Let $M$ be almost invertible. Applying Gram-Schmidt with positive
normalizations on its first $n-1$ columns we obtain the first $n-1$ columns of
both $Q$ and $R$, as well as those of $Q_\star$ and $R_\star$. The last column
$v = Q_\star e_n$ of $Q_\star$ is already well defined, by orthonormality and
the fact that $\det Q_\star = 1$. Now, set $R_\star = M (Q_\star)^\ast$. The
positivity of $R_{n,n}$ specifies whether the last column of $Q$ is $v$ or
$-v$. Smoothness
is clear by construction.

If $M$ is invertible, $\det M = \det Q_\star  \det R_\star$ implies that the
last diagonal entry of $R_\star$ has the same sign of $\det M$: the relations
between the factorizations then follow. If $M$ is not invertible, the relation
among determinants implies $(R_\star)_{n,n}=0$.
\qed

If all subdiagonal entries of $T$ are nonzero, $T$ is
an \textit{unreduced} matrix; otherwise, $T$ is {\it reduced}.
Notice that an unreduced tridiagonal matrix is almost invertible:
indeed, the block formed by
rows $2,\ldots,n$ and columns $1,\ldots,n-1$ is a an upper triangular matrix
with nonzero diagonal entries, and therefore, invertible.

We consider the {\it shifted $QR$ step} and its {\it signed} counterpart,
\[ \Phi(T, s) = Q^\ast T Q, \quad \Phi_\star(T, s) = Q_\star^\ast T Q_\star, \]
where $T - s I = QR$ and $T - s I = Q_\star R_\star$.
Let $\dom(\Phi)$ be the set of pairs  $(T,s) \in \Tt \times \RR$
for which $T- sI$ is invertible:
$\dom(\Phi)$ is open and dense in $\Tt \times \RR$
and, from the Gram-Schmidt algorithm, $\Phi$ is smooth in $\dom(\Phi)$.
Similarly, the above proof shows that $\Phi_\star$ is smooth
in $\dom(\Phi_\star)$,
with $(T,s)  \in \dom(\Phi_\star)$ if $T- sI$ is  almost invertible.
Clearly, $\dom(\Phi)$ is strictly contained in $\dom(\Phi_\star)$.

\begin{lemma}
\label{lemma:basic}
For $(T,s) \in \dom(\Phi)$ (resp.  $\dom(\Phi_\star)$),
we have $\Phi(T,s) \in \Tt$ (resp. $\Phi_\star(T,s) \in \Tt$).
The spectra of $T$, $\Phi(T,s)$ and $\Phi_\star(T,s)$ are equal.
In the appropriate domains,
for $T - sI = QR = Q_\star R_\star$ and $i = 1, 2, \ldots, n-1$,
\[ (\Phi(T,s))_{i+1,i} =
\frac{(R)_{i+1,i+1}}{(R)_{i,i}} \; (T)_{i+1,i}, \quad
(\Phi_\star(T,s))_{i+1,i} =
\frac{(R_\star)_{i+1,i+1}}{(R_\star)_{i,i}} \; (T)_{i+1,i}. \]
Thus, the top $n-2$ subdiagonal entries of $T$,
$\Phi(T,s)$ and $\Phi_\star(T,s)$ have the same sign;
also, $\sign\, (T)_{n,n-1}  = \sign\, (\Phi(T,s))_{n,n-1}$.

\end{lemma}

\smallskip

\proof
We prove the statements for $\Phi_\star$; the others are then easy.

For a pair $(T,s) \in \dom (\Phi) \subset \dom (\Phi_\star)$,
there are two expressions for $\Phi_\star(T,s)$:
\[ \Phi_\star(T,s) = Q_\star^\ast T Q_\star = R_\star T R_\star^{-1},
\quad \hbox{where}\quad T - sI = Q_\star R_\star. \]
From the first equality, $\Phi_\star(T,s)$ is symmetric and
from the second, $\Phi_\star(T,s)$ is an upper Hessenberg matrix
so that $\Phi_\star(T,s) \in \Tt$ is similar to $T$.
More generally, for $(T,s) \in \dom (\Phi_\star)$
we still have
\[ \Phi_\star(T,s) = Q^\ast_\star T Q_\star , \quad
\Phi_\star(T,s) R_\star = R_\star T \]
and therefore $\Phi_\star(T,s) \in \Tt$ is similar to $T$.
Compute the $(i+1,i)$ entry of the second equation above to obtain
$(\Phi_\star(T,s))_{i+1,i} \; (R_\star)_{i,i} =
(R_\star)_{i+1,i+1} \; (T)_{i+1,i}$,
completing the proof.
\qed

The following result describes the behavior of $\Phi_\star$
at points not in $\dom(\Phi)$,
which will play an important role throughout the paper.

\begin{lemma}
\label{lemma:sislambda}
If $(T,s) \in \dom(\Phi_\star) \smallsetminus \dom(\Phi)$ then
\[ \bbb(\Phi_\star(T,s)) = (\Phi_\star(T,s))_{n,n-1} = 0, \quad
(\Phi_\star(T,s))_{n,n} = s. \]
At a point $(T,s) \in \dom(\Phi_\star)$ with
$\bbb(T) = 0$ and $s = (T)_{n,n}$ we
have $\grad(\bbb \circ \Phi_\star) = 0$.
\end{lemma}

\proof
Since $T - sI = Q_\star R_\star = R_\star^\ast Q_\star^\ast$ is not invertible
then $(R_\star)_{n,n} = 0$ and therefore $R_\star^\ast e_n = 0$.
Thus $v = (Q_\star^\ast)^{-1} e_n = Q e_n$ satisfies $(T - sI)v = 0$.
We then have
$\Phi_\star(T,s) e_n = Q^\ast T Q e_n = Q^\ast T v = Q^\ast (sv) = se_n$,
proving the first claim.
For the second claim, since $T - sI$ is almost invertible,
$(R_\star)_{i,i} > 0$ for $i < n$.
From the previous lemma,
\[ (\bbb \circ \Phi_\star)(T,s) =
\frac{(R_\star)_{n,n}}{(R_\star)_{n-1,n-1}} \; \bbb(T)\,; \]
if $\bbb(T) = 0$ and $s = (T)_{n,n}$ then
$(R_\star)_{n,n} = 0$ and
$\bbb \circ \Phi_\star$ is a product of two smooth functions,
both zero, yielding $\grad(\bbb \circ \Phi_\star) = 0$.
\qed

The operation of changing subdiagonal signs, i.e.,
of conjugation by some $E \in \Ee$,
behaves well with respect to $\Phi$ and $\Phi_\star$.
For $1 \le j < n$, let $E_j \in \cE$ be defined by
\[ (E_j)_{i,i} = \begin{cases} 1, & i \le j, \\ -1, & i > j. \end{cases} \]
Together with $-I$, the matrices $E_i$ generate $\cE$.
For $T \in \cT$, $\eta_i(T) = E_i T E_i$ differs from $T$
only in the sign of the $i$-th subdiagonal coordinate:
$(\eta_i(T))_{i+1,i} = - (T)_{i+1,i}$.
The nontrivial involutions in $\ILa$ are therefore generated
by $\eta_i$, $1 \le i < n$.

\begin{lemma}
\label{lemma:ETE}
The domains $\dom(\Phi)$ and $\dom(\Phi_\star)$
are $\cE$-invariant and
\[ \Phi(\eta(T),s)= \eta(\Phi(T,s)), \quad
\Phi_\star(\eta(T),s)= \eta(\Phi_\star(T,s)). \]
If $\det (T - s I) > 0 $ then $\Phi(T,s)=\Phi_\star(T,s)$;
if  $\det (T- s I) < 0 $, $\Phi(T,s)= \eta_{n-1}(\Phi_\star(T,s))$;
if  $\det (T- s I) = 0 $ and $(T,s) \in \dom(\Phi_\star)$,
then $\bbb(\Phi_\star(T,s))=0$.
\end{lemma}


\proof
For $(T,s) \in \dom(\Phi)$,
the matrices $T - sI$ and $E(T - sI)E$ are both invertible.
The $QR$ factorization $T -sI = Q R$ yields
$ETE - E(sI)E = (EQE)(ERE)$,
preserving the positivity of the diagonal entries of the triangular part, so
\[ \Phi(ETE,s) = (EQE)^\ast ETE (EQE) = E Q^\ast T Q E = E \Phi(T,s) E.\]
The argument is similar for $\Phi_\star$.
The claims for $T-sI$ invertible follow from the relation
between $Q$ and $Q_\star$ in Proposition \ref{prop:existence};
the case $\det(T - sI) = 0$ is a repetition of Lemma \ref{lemma:sislambda}.
\qed

We are only interested in the case when the spectrum of $T$ is simple,
since a double eigenvalue implies reducibility.
Since either version of shifted $QR$ step preserves spectrum,
restriction defines smooth
maps
$\Phi: (\ILa \times \RR) \cap \dom(\Phi) \to \ILa$ and
$\Phi_\star: (\ILa \times \RR) \cap \dom(\Phi_\star) \to \ILa$.

Still in $\ILa$, it is convenient to consider
the \textit{step} $F_s(T) = \Phi_\star(T,s)$.
For $s$ not an eigenvalue of $\Lambda$, the domain of $F_s$ is $\ILa$.
The natural domain for $F_{\lambda_i}$ instead
is the \textit{deflation domain} $\Dd_\Lambda^i$,
the open dense subset of $\ILa$ of matrices $T$
for which $T - \lambda_i I$ is almost invertible.
In other words, $T \in \Dd_\Lambda^i$ if and only if
$\lambda_i$ is an eigenvalue of the lowest irreducible block of $T$.

The definition of the step $F_s$ differs from the usual one
in that we use $\Phi_\star$ instead of $\Phi$.
Given Lemma \ref{lemma:ETE},
considerations about deflation are unaffected
and our choice has the advantage of being smooth
(and well defined) in $\Dd_{\Lambda}^i$.

The ($i$-th) \textit{deflation set} is
\[ \Dd_{\Lambda,0}^i =
\left\{ T \in \ILa\;|\; \bbb(T)=0, (T)_{n,n} = \lambda_i \right\}. \]
Since the spectrum of $\Lambda$ is simple,
$\Dd_{\Lambda,0}^i \subset \Dd_\Lambda^i$.
Also, if $i \ne j$ then $\Dd_\Lambda^i \cap \Dd_{\Lambda,0}^j = \emptyset$.
We saw in Lemma \ref{lemma:sislambda} that
when the shift is taken to be an eigenvalue, a single step deflates a matrix,
i.e., that the image of $F_{\lambda_i}$ is contained in $\Dd^i_{\Lambda,0}$:
we shall see in Proposition \ref{prop:Fs} that this image
is in fact equal to $\Dd^i_{\Lambda,0}$.

\begin{prop}
\label{prop:Fs}
If $s$ is not an eigenvalue of $\Lambda$,
the map $F_s: \ILa \to \ILa$ is a diffeomorphism.
The image of $F_{\lambda_i}: \Dd_\Lambda^i \to \ILa$ is $\Dd_{\Lambda,0}^i$.
The restriction
$F_{\lambda_i}|_{\Dd^i_{\Lambda,0}}: \Dd^i_{\Lambda,0} \to \Dd^i_{\Lambda,0}$
is a diffeomorphism.
\end{prop}

\proof
If $s$ is not an eigenvalue,
compute $F_s^{-1}(T)$ by factoring $T - sI$ as $RQ$,
$R$ upper triangular with the first $n-1$ diagonal entries positive
and $Q \in SO(n)$:
we claim that $F_s(T_0) = T$ for $T_0 = QR + sI$, proving
that $F_s$ is a diffeomorphism.
Indeed, $QR = T_0 - sI$ is a $Q_\star R_\star$
factorization and thus $F_s(T_0) = Q^\ast T_0 Q = T$.

From the last sentence of Section 2,
the image of $F_{\lambda_i}$ is contained in
$\Dd^i_{\Lambda,0} \subset \Dd^i_{\Lambda}$.
The fact that the restriction of $F_{\lambda_i}$ to $\Dd^i_{\Lambda,0}$ is
a diffeomorphism is proved as in the previous paragraph.
\qed

Commutativity of steps is well known
and related to the complete integrability
of the interpolating Toda flows
(\cite{Flaschka}, \cite{LT}, \cite{Moser}, \cite{Parlett}).
For the reader's convenience we provide a quick proof.

\begin{prop}
\label{prop:commute}
Steps commute:
$F_{s_0} \circ F_{s_1} = F_{s_1} \circ F_{s_0}$
in the appropriate domains.
\end{prop}

The domain of $F_{s_0} \circ F_{s_1} = F_{s_1} \circ F_{s_0}$ is $\ILa$ if
neither $s_0$ nor $s_1$ is an eigenvalue,
$\Dd^i_\Lambda$ if $s_0 = \lambda_i$
and $s_1$ is not an eigenvalue (or vice-versa) and
the empty set in the rather pointless case
$s_0 = \lambda_i$, $s_1 = \lambda_j$, $i \ne j$.

\smallskip

\proof
We prove commutativity only when $s_0$ and $s_1$ are not eigenvalues;
the other cases follow easily.
Consider $Q_\star R_\star$ factorizations
\begin{gather*}
T - s_0 I = Q_0 R_0, \quad T - s_1 I = Q_1 R_1, \\
(T - s_0 I)(T - s_1 I) = (T - s_1 I)(T - s_0 I) = Q_2 R_2.
\end{gather*}
For $F_{s_0}(T) - s_1 = Q_0^\ast(T - s_1)Q_0 = Q_3 R_3$,
we have $F_{s_1}(F_{s_0}(T)) = Q_3^\ast F_{s_0}(T) Q_3 =
Q_3^\ast Q_0^\ast T Q_0 Q_3$.
Thus
\[ Q_0^\ast(T - s_1)Q_0R_0 = Q_0^\ast (T - s_1 I)(T - s_0 I) =
Q_0^\ast Q_2 R_2 = Q_3 R_3 R_0 \]
and therefore $Q_0^\ast Q_2 = Q_3$ and
$F_{s_0}(F_{s_1}(T)) = Q_2^\ast T Q_2$.
\qed

\section{Simple shift strategies}

The point of using a shift strategy is to accelerate deflation,
ideally by choosing $s$ near an eigenvalue of $T$.
A \textit{simple shift strategy} is an $\cE$-invariant function
$\sigma: \ILa \to \RR$ such that
there exists $C_\sigma > 0$ such that
for all $T \in \ILa$ there is an eigenvalue $\lambda_i$ with
$| \sigma(T) - \lambda_i | \le C_\sigma |\bbb(T)|$.
In particular, if $T \in \Dd_{\Lambda,0}^i$ then $\sigma(T) = \lambda_i$.

The \textit{step} associated with a (simple) shift strategy $\sigma$
is $F_\sigma$, defined by $F_\sigma(T) = F_{\sigma(T)}(T)$.
The natural domain for $F_\sigma$ is the set of matrices $T$
for which $T - \sigma(T) I$ is almost invertible.
From Section 2, it includes all unreduced matrices
and open neighborhoods of each deflation set $\Dd^i_{\Lambda,0}$.
We shall also see in Section 6 that it contains
a dense open subset $\Uu_{\Lambda,\epsilon}$ of $\ILa$ invariant under $F_\sigma$.  
A more careful description of this domain will not be needed.


Quoting Parlett (\cite{Parlett}), there are shifts for all seasons.
Let $\rho$ be \textit{Rayleigh's shift}:
$\rho(T) = (T)_{n,n}$.
Denote the bottom $2 \times 2$ diagonal principal minor
of a matrix $T \in \Tt$ by $\hat{T}$:
\textit{Wilkinson's shift} $\omega(T)$ is the eigenvalue
of $\hat T$ closer to $(T)_{n,n}$
(in case of draw, take the smallest eigenvalue).

\begin{lemma}
\label{lemma:2sqrt2}
The functions $\rho$ and $\omega$ is are simple shift
strategies with $C_\rho = \sqrt{2}$ and $C_\omega = 2\sqrt{2}$.
\end{lemma}

We use here the Wielandt-Hoffman theorem
(for a simple proof using the Toda dynamics, see \cite{DRTW}):
if $S, T \in \Tt$ have eigenvalues
$\sigma_i$ and $\lambda_i$ in increasing order then
\[ \sum_i |\sigma_i - \lambda_i|^2 \le \trace((S-T)^2). \]

\proof
Invariance is trivial for $\rho$;
for $\omega$, it follows from the fact that
changing signs of off-diagonal
entries of a $2 \times 2$ matrix does not change its spectrum.

Let $B = e_n e_{n-1}^\ast + e_{n-1} e_n^\ast$
and $S = T - \bbb(T) B$ so that $\rho(T) = (T)_{n,n}$ is an eigenvalue of $S$.
From the Wielandt-Hoffman theorem, for some $i$,
\[ |\rho(T) - \lambda_i| \le \sqrt{2} \; \left|\bbb(T)\right|, \]
proving that $C_\rho = \sqrt{2}$.
Apply again the Wielandt-Hoffman theorem
to the $2 \times 2$ trailing principal minors of $S$ and $T$
to deduce that
\[ |(T)_{n,n} - \omega(T)| \le \sqrt{2} \; \left|\bbb(T)\right|. \]
We thus have $| \omega(T) - \lambda_i| \le 2\sqrt{2} \; \left|\bbb(T)\right|$
and $C_\omega = 2\sqrt{2}$, as desired.
\qed

Another example of (simple) shift strategy,
the \textit{mixed Wilkinson-Rayleigh strategy},
uses Wilkinson's shift unless the matrix is already near deflation,
in which case we use Rayleigh's:
\[ \sigma(T) = \begin{cases}
\rho(T),&|(T)_{n,n-1}| < \epsilon, \\
\omega(T),&|(T)_{n,n-1}| \ge \epsilon;
\end{cases} \]
here $\epsilon > 0$ is a small constant.

Simple shift strategies are not required to be continuous and
$\omega$ is definitely not.
For a simple shift strategy $\sigma$,
let $\calS_\sigma \subset \ILa$ be the \textit{singular support} of $\sigma$,
i.e., a minimal closed set on whose complement $\sigma$ is smooth.
For example, $\calS_\omega$ is the set of matrices $T \in \ILa$
for which the two eigenvalues $\omega_-(T)$ and $\omega_+(T)$ of $\hat T$
are equidistant from $(T)_{n,n}$, or,
equivalently, for which $(T)_{n,n} = (T)_{n-1,n-1}$.
The set $\calS_\sigma$ will play an important role later.

We consider the phase portrait of $F_\omega$ for $3 \times 3$ matrices.
In this case,
the reader may check that the domain of $F_\omega$ is the full set $\ILa$.
Let $\JLao \subset \ILa$ be set of \textit{Jacobi matrices}
similar to $\Lambda$, i.e., matrices $T \in \ILa$
with strictly positive subdiagonal entries.
Recall that the closure $\JLa \subset \ILa$ is diffeomorphic to a hexagon,
the permutohedron in this dimension.
The set $\JLa$ is not invariant under $F_\omega$
but we may define $\tilde F_\omega(T)$ with $\tilde F_\omega: \JLa \to \JLa$
by dropping signs of subdiagonal entries of $F_\omega(T)$.
As discussed above, this standard procedure is mostly harmless.

\begin{figure}[ht]
\begin{center}
\psfrag{Yy}{$\calS_\omega$}
\epsfig{height=28mm,file=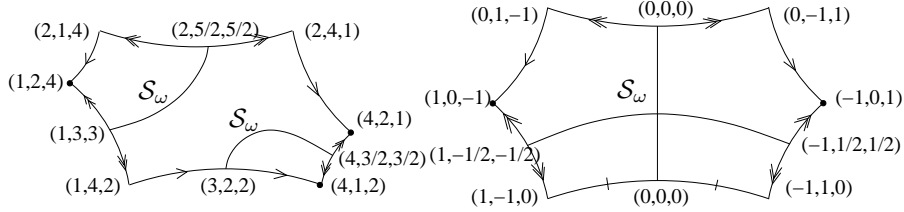}
\end{center}
\caption{\capsize The phase space of Wilkinson's step for $n=3$.}
\label{fig:wil}
\end{figure}

Two examples of $\tilde F_\omega$ are given in Figure \ref{fig:wil},
which represent $\JLa$ for the $\Lambda= \diag(1,2,4)$ on the left
and $\Lambda=\diag(-1,0,1)$ on the right.
The vertices are the six diagonal matrices
similar to $\Lambda$ and the edges consist of reduced matrices.
Labels indicate the diagonal entries of the corresponding matrices.
Three edges form $\Dd^i_{\Lambda,0} \cap \JLa$:
they alternate, starting from the bottom horizontal edge on both hexagons.
The set $\calS \cap \JLa$ is indicated in both cases.


Vertices are fixed points of $\tilde F_\omega$ and boundary edges are invariant sets.
A simple arrow indicates the motion of the points $\tilde F_\omega^k(T)$ along the edge.
Points $T$ on an arc with a double arrow are taken to a diagonal matrix in a single step:
the arc points to $\tilde F_\omega(T)$.
Arcs marked with a transversal segment consist of fixed points of $\tilde F_\omega$.

Points on both sides of $\calS_\omega$ are taken far apart:
there is a jump discontinuity along $\calS_\omega$.
From Theorem \ref{theo:squeeze},
the decay of the bottom subdiagonal entry under Wilkinson's step
away from $\calS_\omega \cap \Dd_{\Lambda,0}$ is cubic.
As discussed in \cite{LST2}, near $\calS_\omega \cap \Dd_{\Lambda,0}$
this decay is quadratic, but not cubic.
For the left hexagon, cubic convergence occurs in the long run
because the sequence $\tilde F_\omega^k(T)$ stays close to this intersection
only for a few values of $k$,
illustrating Theorem \ref{theo:big}.

In the case $\Lambda= \diag(-1,0,1)$, the bottom edge consists of fixed points.
As mentioned in the Introduction (see Figure \ref{fig:xaxbxc}),
this case has a special asymptotic behavior:
the (fixed) point labeled by $(0,0,0)$ is
the central point of the set $\cX$.
If $T \in \cX$ then the sequence $(\tilde F_\omega^k(T))$
is contained in $\cX$ and converges to the central point
at a strictly quadratic rate.


\section{Tubular coordinates}

Recall that a map $\Pi: X \to Y \subset X$ is a \textit{projection}
if $\Pi(X) = Y$ and $\Pi \circ \Pi = \Pi$.
Instead of using abstract topological facts to prove
the existence of some projection $\Dd_\Lambda^i \to \Dd^i_{\Lambda,0}$
we prefer to construct a specific projection which
works well with the $QR$ steps.
The map $F_{\lambda_i}: \Dd_\Lambda^i \to \Dd^i_{\Lambda,0}$ 
is not a projection but can, using Proposition \ref{prop:Fs},
be used to define one: the \textit{canonical projection} 
$\Pi_i: \Dd_\Lambda^i \to \Dd^i_{\Lambda,0}$,
\[ \Pi_i(T) =
(F_{\lambda_i}|_{\Dd^i_{\Lambda,0}})^{-1}(F_{\lambda_i}(T)). \]

\begin{prop}
\label{prop:pia}
The map $\Pi_i$ is a
smooth projection which commutes with steps:
$\Pi_i(F_s(T)) = F_s(\Pi_i(T))$ provided
$s$ is not an eigenvalue of $\Lambda$ different from $\lambda_i$.
\end{prop}

\proof
The map $\Pi_i$ is clearly smooth and, for $T \in
\Dd^i_{\Lambda,0}$, we have \[ \Pi_i(T) =
(F_{\lambda_i}|_{\Dd^i_{\Lambda,0}})^{-1}(F_{\lambda_i}(T)) = T, \] proving
that $\Pi_i$ is a projection.  Commutativity follows from
Proposition \ref{prop:commute}.
\qed

For a diagonal matrix $\Lambda$ with simple spectrum and $\epsilon > 0$,
the \textit{deflation neighborhood}
$\Dd_{\Lambda,\epsilon}\subset \ILa$
is the closed set of matrices $T \in \ILa$ with $|\bbb(T)| \le \epsilon$.
This notation is consistent with $\Dd_{\Lambda,0}$
for the deflation set.
As we shall see in Propositions \ref{prop:deflation} and \ref{prop:epsilonF},
for sufficiently small $\epsilon > 0$ the set $\Dd_{\Lambda,\epsilon}$
has connected components $\Dd^i_{\Lambda,\epsilon} \subset \Dd^i_\Lambda$,
$\Dd^i_{\Lambda,\epsilon} \supset \Dd^i_{\Lambda,0}$,
which are invariant under steps $F_s$ for shifts $s$ near $\lambda_i$, i.e.,
$F_s(\Dd^i_{\Lambda,\epsilon}) \subset \Dd^i_{\Lambda,\epsilon}$.
The sets $\Dd^i_{\Lambda,\epsilon}$ are therefore
also invariant under $F_\sigma$.

Denote the distance between a matrix $T$
and a compact set of matrices $\cN$ by
$\dist(T,\cN) = \min_{S \in \cN} \| T - S \|$.
Let $\gamma = \min_{i \ne j} |\lambda_i - \lambda_j|$ be
the {\it spectral gap} of $\Lambda$
and $B = e_n e_{n-1}^\ast + e_{n-1} e_n^\ast$.

Recall that if $\cN$ is a submanifold of codimension $k$ of $\cM$
then a \textit{closed tubular neighborhood} of $\cN$
consists of a closed neighborhood $\cN_\epsilon$ of $\cN$
and a diffeomorphism
$\zeta: \cN_\epsilon \to \cN \times \BB^k_\epsilon$ with
$\zeta(x) = (x,0)$ for $x \in \cN$
(here $\BB^k_\epsilon \subset \RR^k$ is the closed ball
of radius $\epsilon$ around the origin).
Given $x \in \cN$, the preimage
$\zeta^{-1}(\{x\} \times \BB^k_\epsilon)$
is a manifold with boundary of dimension $k$,
the \textit{fiber} through $x$.
We now construct tubular neighborhoods
of the deflation sets $\Dd^i_{\Lambda,0}$;
here the codimension is $k = 1$.

\begin{prop}
\label{prop:deflation}
Each $\Dd_{\Lambda,0}^i \subset \ILa$ is a
compact submanifold of codimension $1$ diffeomorphic to $\ILai$,
where $\Lambda_i = \diag(\lambda_1, \ldots, \lambda_{i-1},
\lambda_{i+1}, \ldots \lambda_n)$.
There exists $\epsilon_\tub > 0$ such that
for $\epsilon \in (0,\epsilon_\tub)$:
\begin{enumerate}[(a)]
\item{the connected components
$\Dd_{\Lambda,\epsilon}^i$ of $\Dd_{\Lambda,\epsilon}$ consist of matrices
$T \in \Dd_{\Lambda,\epsilon}$ for which
$|(T)_{n,n} - \lambda_i| < \sqrt{2} \;\epsilon$;}
\item{the map
$\zeta: \Dd^i_{\Lambda,\epsilon} \to
\Dd^i_{\Lambda,0} \times [-\epsilon, \epsilon]$
given by $\zeta(T) = (\Pi_i(T), \bbb(T))$
is a closed tubular neighborhood of $\Dd_{\Lambda,0}$;}
\item{there is a constant $C_\bbb>0$ such that for all
$T \in \Dd^i_{\Lambda,\epsilon}$,
\[ | b(T)| \le \dist(T,\Dd^i_{\Lambda,0}) \le
\| T - \Pi_i(T)\| \le C_\bbb |b(T)|.\]}
\end{enumerate}
\end{prop}

\smallskip
\proof
We first show that the gradient of the restriction
$\bbb|_{\ILa}$ at a point $T_{\Dd} \in \Dd_{\Lambda,0}$ is not zero.
Consider the characteristic polynomial
along the line $T_{\Dd} + tB$:
this is a smooth even function of $t$
and therefore $B$ is tangent to $\ILa$ at $T_{\Dd}$,
the point on which $t=0$.
On the other hand, the directional derivative of $\bbb$
along  the same line equals 1.
Thus $\Dd_{\Lambda,0} \subset \ILa$ is a
submanifold of codimension $1$.
The diffeomorphism with $\ILai$ takes $T$ to
$\hat T$, the leading $(n-1) \times (n-1)$ principal minor of $T$.

Assume $\epsilon < \gamma/(2 \sqrt{2})$.
Consider matrices $T \in \Dd_{\Lambda,\epsilon}$ and $S =T - \bbb(T) B$,
so that $(T)_{n,n}$ is an eigenvalue of $S$.
By the Wielandt-Hoffman theorem,
there exists an index $i$ for which
$|(T)_{n,n} - \lambda_i| < \sqrt{2} \epsilon$,
defining the sets
$\Dd^i_{\Lambda,\epsilon}$
(at this point we do not yet know
that $\Dd^i_{\Lambda,\epsilon}$ is connected).



For $T_{\Dd} \in \Dd^i_{\Lambda,0}$,
the derivative $D\Pi_i(T_{\Dd})$
equals the identity on the subspace tangent
to $\Dd^i_{\Lambda,0}$ and has a kernel of dimension $1$.
Thus, for sufficiently small $\epsilon_\tub$, item (b) holds.
This also proves that each $\Dd^i_{\Lambda,\epsilon}$ is connected,
completing the proof of item (a).

The first two inequalities in (c) are trivial. Now
\[ \| T - \Pi_i(T) \| =
\| \zeta^{-1}(\Pi_i(T),\bbb(T)) -
\zeta^{-1}(\Pi_i(T),0) \| \le
C_\bbb | \bbb(T) |, \]
where the derivative of $\zeta^{-1}(T_{\Dd}, \delta)$
with respect to the second
coordinate is bounded by $C_\bbb$ on the compact set
$\Dd_{\Lambda,0} \times [-\epsilon_\tub, \epsilon_\tub]$.
\qed


The diffeomeorphism $\zeta$ defines \textit{tubular coordinates}
for $T \in \Dd^i_{\Lambda,\epsilon}$:
the matrix $\Pi_i(T) \in \Dd^i_{\Lambda,0} \approx \Tt_{\Lambda_i}$
and $\bbb(T)$.
Under tubular coordinates, $QR$ steps with shift are given by a simple formula.

\begin{coro}
\label{coro:tubularQR}
Consider $\Lambda$, $i$ and $\epsilon \in (0,\epsilon_\tub)$. Then
\begin{align*}
\zeta \circ F_s \circ \zeta^{-1}:
\Dd^i_{\Lambda,0} \times [-\epsilon, \epsilon] &\to
\Dd^i_{\Lambda,0} \times [-\epsilon, \epsilon] \\
(T,b) &\mapsto \left( F_s(T), \frac{(R_\star)_{n,n}}{(R_\star)_{n-1,n-1}} \,b \right)
\end{align*}
where $\zeta^{-1}(T,b) - sI = Q_\star R_\star$.
\end{coro}

\proof
This follows directly from Lemma \ref{lemma:basic} and Propositions \ref{prop:pia} and \ref{prop:deflation}.
\qed

\section{Convergence to deflation}

Sufficiently thin deflation neighborhoods $\Dd_{\Lambda,\epsilon}^i$ are
invariant under $F_s$ for $s \approx \lambda_i$.

\begin{prop}
\label{prop:epsilonF}
Given $C > 0$, there exists
$\epsilon_\inv \in (0, \epsilon_\tub)$,
such that for any $\epsilon \in (0,\epsilon_\inv)$ and
$s \in [\lambda_i - C\,\epsilon,\lambda_i + C\,\epsilon]$ we have
$F_s(\Dd^i_{\Lambda,\epsilon}) \subset \interior(\Dd^i_{\Lambda,\epsilon/2})$.

For a simple shift strategy
$\sigma: \ILa \to \RR$,
there exists $\epsilon_\inv > 0$ such that
if $\epsilon \in (0,\epsilon_\inv)$ then
$F_\sigma(\Dd^i_{\Lambda,\epsilon}) \subset
\interior(\Dd^i_{\Lambda,\epsilon/2})$.
\end{prop}

In particular, $F_s$ is well defined in $\Dd^i_{\Lambda,\epsilon}$ for
$\epsilon \in (0,\epsilon_\inv)$.

\smallskip

\proof
Recall that $F_s(\Dd_{\Lambda,0}^i) = \Dd_{\Lambda,0}^i$.  
From Lemma \ref{lemma:sislambda}, the derivative of $\bbb \circ \Phi_\star$ is zero
at $\Dd_{\Lambda,0}^i \times \{ \lambda_i \}$.
Compactness of $\Dd_{\Lambda,0}^i$ thus implies that
in a sufficiently small neighborhood of
$\Dd_{\Lambda,0}^i \times \{ \lambda_i \}$
we have $|\bbb(F_s(T))| \le |\bbb(T)|/3$.

Now consider a simple shift strategy $\sigma$:
there exists $C_\sigma > 0$ such that
$|\sigma(T) - \lambda_i| < C_\sigma \bbb(T)$;
apply the first statement with $C = C_\sigma$.
\qed


Thus, $F_\sigma$ squeezes neighborhoods
$\Dd^i_{\Lambda,\epsilon}$ at least linearly.
Equivariance and smoothness imply an estimate stronger
than that in the definition of simple shift strategy.
We do not want to assume, however,
that $\Dd_{\Lambda,0} \cap \calS_{\sigma} = \emptyset$:
after all, this is not true even for Wilkinson's shift.
We need a more careful statement.

\begin{lemma}
\label{lemma:wlip}
Consider a shift strategy $\sigma$ and
$\epsilon_\inv$ as in Proposition \ref{prop:epsilonF}.
For a compact set $\Kk \subset
\Dd^i_{\Lambda,\epsilon_\inv} \smallsetminus
(\Dd^i_{\Lambda,0} \cap \calS_\sigma)$,
there exists $C_{\Kk}$ such that for all $T \in \Kk$ we have
$|\sigma(T) - \lambda_i| \le C_{\Kk} \bbb(T)^2 $.
\end{lemma}

\proof
Let $\Kk_{\Dd} = \Kk \cap \Dd^i_{\Lambda,0}$;
enlarge $\Kk_{\Dd}$ along $\Dd_{\Lambda,0}^i$
to obtain another compact set
$\Kk_{1} \subset \Dd_{\Lambda,0}^i \smallsetminus \calS_\sigma$,
$\Kk_{\Dd} \subset \interior_{\Dd_{\Lambda,0}^i}(\Kk_{1})$.
Fatten $\Kk_1$ along fibers to define
$\tilde\Kk_1 = \zeta^{-1}(\Kk_1 \times [-\epsilon,\epsilon])$,
$\epsilon \in (0,\epsilon_\inv)$, which, without loss,
still avoids $\calS_\sigma$.
For each $T_\Dd \in \Kk_1$,
consider the function $h_{T_{\Dd}}(b) = \sigma(\zeta^{-1}(T_\Dd,b))$,
obtained by restricting $\sigma$ to a fiber of $\Dd_{\Lambda,\epsilon}^i$.
Each $h_{T_{\Dd}}$ is smooth and even
and therefore satisfies
$|h_{T_{\Dd}}(b) - \lambda_i| \le C_{T_{\Dd}} |b|^2$.
By compactness, there exists $C_{\Kk_1}$ such that
$|h_{T_{\Dd}}(b) - \lambda_i| \le C_{\Kk_1} |b|^2$ for all $T_{\Dd} \in \Kk_1$.
In other words, there exists $C_{\tilde\Kk_1}$ such that
$|\sigma(T) - \lambda_i| \le C_{\tilde\Kk_1} |\bbb(T)|^2$
for all $T \in \tilde\Kk_1$.
The estimate for $T \notin \tilde\Kk_1$ is trivial.
\qed


\smallskip

{\noindent \bf Proof of Theorem \ref{theo:squeeze}:}
Take $\epsilon = \epsilon_\inv$ as in Proposition \ref{prop:epsilonF}
so that $\Dd^i_{\Lambda,\epsilon}$ is invariant under $F_\sigma$.

Let $\varphi = \bbb \circ \Phi_\star$.
We compute the Taylor expansion of $\varphi(T,s)$ at $(T_{\Dd},\lambda_i)$,
$T_{\Dd} \in \Dd_{\Lambda,0}^i$:
from Lemma \ref{lemma:sislambda}, the gradient of $\varphi$
at $(T_{\Dd},\lambda_i)$ is zero.
Thus, up to a third order remainder,
\begin{align*}
\varphi(T, s) &=
\varphi(T_{\Dd}, \lambda_i) +
\frac{1}{2}\varphi_{T,T}(T_{\Dd}, \lambda_i)(T - T_{\Dd},T - T_{\Dd}) + \\
&\quad + \varphi_{T,s}(T_{\Dd}, \lambda_i)(T - T_{\Dd},s - \lambda_i) +
\frac{1}{2} \varphi_{s,s}(T_{\Dd}, \lambda_i)(s - \lambda_i,s - \lambda_i)
+ \\ &\quad + \hbox{Rem}_3 ( T - T_{\Dd},s - \lambda_i).
\end{align*}
Now, $\varphi(T_{\Dd}, \lambda_i)=0$
and, again from Lemma \ref{lemma:sislambda},
$\varphi(T,\lambda_i) = 0$ for all $T \in \ILa$,
hence $\varphi_{T,T}(T_{\Dd}, \lambda_i)=0$.
Let $C_\sigma$ be the constant in
the definition of a simple shift strategy.
By compactness, there exists $C_1 >0$ such that
for all $T_{\Dd} \in \Dd^i_{\Lambda,0}$,
$T \in \Dd^i_{\Lambda,\epsilon}$ and
$s \in [\lambda_i - C_\sigma\,\epsilon, \lambda_i + C_\sigma\,\epsilon]$,
we have
\[ | \varphi(T,s) | \le
C_1 |s - \lambda_i| ( \| T - T_{\Dd} \| + |s - \lambda_i|) \]
We now apply this estimate for $T_{\Dd} = \Pi_i(T)$,
where $T \in \Dd^i_{\Lambda,\epsilon}$.
By Proposition \ref{prop:deflation}, since $\epsilon < \epsilon_\tub$,
$ \| T - T_{\Dd} \| = \| T - \Pi_i(T)\| \le C_\bbb |\bbb(T)|$
and therefore
\[ | \varphi(T,s) | \le
C_1 |s - \lambda_i| ( C_\bbb |\bbb(T)| + |s - \lambda_i|) \]
implying the quadratic estimate
\[ |\bbb(F_\sigma(T))| =
| \varphi(T,\sigma(T)) | \le
C_1 |\sigma(T) - \lambda_i| ( C_\bbb |\bbb(T)| + |\sigma(T) - \lambda_i|) \le
C_q |\bbb(T)|^2. \]
Using Lemma \ref{lemma:wlip} yields the cubic estimate in (c).
\qed

As a corollary, we obtain the well known fact that, near deflation,
the rate of convergence of Rayleigh's
(as well as the mixed Wilkinson-Rayleigh) strategy
has cubic convergence.  
The rate of convergence for Wilkinson's strategy is subtler.


We construct a larger invariant set for $F_\sigma$.
Let $\Uu_\Lambda \subset \ILa$ be the set of unreduced matrices;
for $\epsilon > 0$,
let $\Uu_{\Lambda,\epsilon} =
\Uu_\Lambda \cup \interior(\Dd_{\Lambda,\epsilon})$.
Notice that $\Uu_{\Lambda,\epsilon}$ is open, dense and path-connected.

\begin{lemma}
\label{lemma:ULaeps}
For a shift strategy $\sigma: \ILa \to \RR$,
$\epsilon_\inv$ as in Proposition \ref{prop:epsilonF}
and $\epsilon \in (0,\epsilon_\inv)$,
the open set $\Uu_{\Lambda,\epsilon}$ is invariant under $F_\sigma$.
\end{lemma}

\proof
If $T \in \Uu_\Lambda$ and $\sigma(T)$ is not in the spectrum then
$F_\sigma(T)$ is (well defined and) unreduced.
If $T \in \Uu_\Lambda$ and $\sigma(T) = \lambda_i$ then
$F_\sigma(T) \in \Dd^i_{\Lambda,0} \subset \Uu_{\Lambda,\epsilon}$.
Finally, if $T \in \interior(\Dd^i_{\Lambda,\epsilon})$ then,
by  Proposition \ref{prop:epsilonF},
$F_\sigma(T) \in \interior(\Dd^i_{\Lambda,\epsilon/2}) \subset
\Uu_{\Lambda,\epsilon}$.
\qed

Notice that we do not assume $\sigma$ or $F_\sigma$ to be continuous.
This shows that for $F_\sigma$ defined from a simple shift strategy $\sigma$
the extra hypothesis in Theorem \ref{theo:connect}, item (b),
actually holds: just take $\cK_i = \Dd^i_{\Lambda,\epsilon}$.

A simple shift strategy $\sigma$ is \textit{deflationary}
if for any $T \in \Uu_{\Lambda,\epsilon_\inv}$
there exists $K \in \NN$  such that
$F_\sigma^K(T) \in \Dd_{\Lambda,\epsilon_\inv}$.
It is now a corollary of Theorem \ref{theo:connect}
and Lemma \ref{lemma:ULaeps}
that continuous simple shift strategies are not deflationary.

Rayleigh's strategy is known not to be deflationary.
The following well known estimate (\cite{HP} and \cite{Parlett}, section 8-10)
implies that Wilkinson's strategy is not only deflationary but uniformly so,
in the sense that there exists $K$ with
$F_\omega^K(\Uu_{\Lambda,\epsilon_\inv}) \subset \Dd_{\Lambda,\epsilon_\inv}$.  As a corollary, the mixed
Wilkinson-Rayleigh strategy is also uniformly deflationary provided
$\epsilon > 0$ is sufficiently small.


\begin{fact}
\label{fact:parlett}
For $T \in  \Tt$ and $k \in \NN$,
\[ |\bbb(F_\omega^k(T))|^3 \leq
\frac{|\bbb(T)^2 (T)_{n-1,n-2}|}{(\sqrt{2})^{k-1}}.\]
\end{fact}

In \cite{Parlett}, the result is shown for unreduced matrices;
the case $T \in \Uu_{\Lambda,\epsilon_\inv}$
follows by taking limits.
Notice that for $T \in \ILa$,
the numerator $|\bbb(T)^2 (T)_{n-1,n-2}|$ is uniformly bounded.


\section{Dynamics for a.p.\ free spectra}

From the previous section, cubic convergence may be lost when the orbit
$F_\sigma^k(T)$ passes near the set $\calS_\sigma \cap \Dd_{\Lambda,0}$.
Our next task is to measure when this happens,
by studying the dynamics associated to a shift strategy
in a deflation neighborhood, i.e.,
the iterates of $F_\sigma:\Dd^i_{\Lambda,\epsilon} \to \Dd^i_{\Lambda,\epsilon}$,
$\epsilon \in (0,\epsilon_\inv)$.
Most of what we need can be read in the projection onto $\Dd^i_{\Lambda,0}$,
where $F_\sigma$ coincides with $F_{\lambda_i}$.

A matrix $T \in \Tt$ with simple spectrum is \textit{a.p.\ free}
if no three eigenvalues are in arithmetic progression and a.p.\ otherwise.
Different kinds of spectra lead to different dynamics:
in this section we handle the a.p.\ free case, clearly a generic restriction.
Let $\tilde{T}$ be the leading principal $(n-1) \times (n-1)$ minor of $T$.
The following result is standard.

\begin{prop}
\label{prop:APfreedynamics}
Let $\Lambda \in \Tt$ be an $n \times n$ diagonal a.p.\ free matrix
with spectrum $\lambda_1 < \cdots < \lambda_n$.
For each $i$,
consider $F_{\lambda_i}: \Dd^i_{\Lambda,0} \to \Dd^i_{\Lambda,0}$
as above. For any ${T} \in \Dd^i_{\Lambda,0}$, the sequence
$(F_{\lambda_i}^k(T))$ converges to a diagonal matrix.
\end{prop}

\proof
The map $F_{\lambda_i}$ on $\Dd^i_{\Lambda,0}$ amounts to a $QR$ step
with shift $\lambda_i$ on $\tilde{T}$,
which has eigenvalues $\lambda_j$, $j \ne i$.
The a.p.\ free hypothesis implies that the absolute values of the
eigenvalues of $\tilde T - \lambda_i I$ are distinct.
If $\tilde T$ is unreduced then, as is well known,
the standard $QR$ iteration converges to a diagonal matrix,
with diagonal entries in decreasing order of absolute value.
More generally, if $\tilde T$ is reduced, apply the above result to each
unreduced sub-block.
\qed

We shall use \textit{height functions} for the $QR$ steps $F_s$, $s$ near $\lambda_i$, i.e.,
functions $H_i: \Dd^i_{\Lambda,\epsilon} \to \RR$
with $H_i(F_s(T)) > H_i(T)$ provided $T$ is not diagonal.
Such height functions and related scenarios
have been considered in \cite{BBR}, \cite{DRTW}, \cite{LT} and \cite{Tomei}.

The matrix $W = \diag(w_1,\ldots,w_n)$ is a \textit{weight matrix}
if $w_1 > \cdots > w_n$.
Since $\Lambda$ is a.p.\ free,
there exists $\epsilon_{\textrm{\rm ap}} \in (0,\epsilon_\inv)$ such that
if $s \in \Ii_i = [\lambda_i - \epsilon_{\textrm{\rm ap}},
\lambda_i + \epsilon_{\textrm{\rm ap}}]$ then
the numbers $|\lambda_j - s|$ are distinct and their order does not depend on $s$.

\begin{prop}
\label{prop:height}
Let $\Lambda$ be an a.p.\ free diagonal matrix,
$W$ a weight matrix and $\epsilon_{\textrm{\rm ap}}$ as above.
For $\delta_H > 0$, set $\lolo_i(x) = \log((x-\lambda_i)^2 + \delta_H)$
and let $H_i: \Dd^i_{\Lambda,\epsilon_{\textrm{\rm ap}}} \to \RR$
be defined by $H_i(T) = \trace(W \lolo_i(T))$.
There exists $\delta_H > 0$ such that
\[ \max_{T \in \partial \Dd^i_{\Lambda,\epsilon_{\textrm{\rm ap}}}} H_i(T)
< \min_{T \in \Dd^i_{\Lambda,0}} H_i(T) \]
and, for any $s \in \Ii_i$,
$H_i$ is a height function
for $F_s: \Dd^i_{\Lambda,\epsilon_{\textrm{\rm ap}}}
\to \Dd^i_{\Lambda,\epsilon_{\textrm{\rm ap}}}$.
\end{prop}

Here,
$\lolo_i(T) = X \diag(\lolo_i(\lambda_1), \ldots, \lolo_i(\lambda_n)) X^{-1}$
for $T = X \Lambda X^{-1}$ so that if $p$ is a polynomial and
$\lolo_i(\lambda_j) = p(\lambda_j)$ for $j = 1, \ldots, n$
then $\lolo_i(T) = p(T)$.
The only conditions on $\lolo_i$ which will be used in the proof are
that $|\lambda_j - \lambda_i| < |\lambda_k - \lambda_i|$ implies
$\lolo_i(\lambda_j) < \lolo_i(\lambda_k)$
and that $\lolo_i(\lambda_i)$ is very negative (for small $\delta_H$).

The proof requires some basic facts about $f$-$Q_\star R_\star$ \textit{steps}
(again related to the integrability of the Toda lattice);
these facts will not be used elsewhere.
For a real diagonal matrix $\Lambda$ with simple spectrum,
let $\Oo_\Lambda$ be the set of
all real symmetric matrices similar to $\Lambda$;
it is well known that $\Oo_\Lambda$ is a smooth compact manifold.
The $f$-$Q_\star R_\star$ \textit{step} applied to a matrix
$S \in \OLa$ is the map $F_f: \Aa_{\Lambda,f} \to \OLa$
defined by $F_f(S) = Q_\star^\ast S Q_\star$,
where $Q_\star$ is obtained from the factorization
$f(S) = Q_\star R_\star$ and $S \in \Aa_{\Lambda,f}$
if and only if $f(S)$ is almost invertible.
If $T \in \ILa \cap \Aa_{\Lambda,f}$ then 
$F_f(T) \in \ILa$ (use the same proof as in Lemma \ref{lemma:basic}).
The maps $F_s: \ILa \to \ILa$ defined above
correspond to restrictions of $F_f$ for $f(x) = x - s$.

%

For a continuous function $h: \RR \to \RR$,
if $S \in \Oo_\Lambda$ then
the matrix function $h(S)$ belongs to $\Oo_\Mu$,
where $\Mu = h(\Lambda)$.
With the obvious abuse of notation,
we have a diffeomorphism $h: \Oo_\Lambda \to \Oo_\Mu$
provided $h$ is injective in the spectrum of $\Lambda$.

\begin{lemma}
\label{lemma:ftilf}
For $h$ injective in the spectrum of $\Lambda$,
consider the diffeomorphism $h: \Oo_\Lambda \to \Oo_\Mu$,
where $\Mu = h(\Lambda)$.
Let $f$ and $\tilde f$ be continuous functions
defined in neighborhoods of the spectra
of $\Lambda$ and $\Mu$, respectively,
satisfying $\tilde f(h(\lambda_j)) = f(\lambda_j)$ for each $j$
with $QR$ steps $F_f:  \Oo_\Lambda \to \Oo_\Lambda$
and $F_{\tilde f}: \Oo_\Mu \to \Oo_\Mu$.
Then $h \circ F_f =  F_{\tilde f} \circ h$.
\end{lemma}

\proof
The hypothesis implies that, for $T \in \Oo_\Lambda$,
$f(T) = \tilde f(h(T)) = QR$ and hence
$F_f(T) = Q^\ast T Q$ and $F_{\tilde f}(h(T)) = Q^\ast h(T) Q$.
Thus $h(F_f(T)) = F_{\tilde f}(h(T))$.
\qed

Let $I_r$ be the $n \times n$ truncated identity matrix,
i.e., $(I_r)_{i,i} = 1$ for $i \le r$, other entries being equal to zero.

\begin{lemma}
\label{lemma:increasing}
Let $\Mu$ be a diagonal matrix with simple spectrum
and $\tilde f: \RR \to \RR$ be a function
for which $\mu_i < \mu_j$ implies
$|\tilde f(\mu_i)| < |\tilde f(\mu_j)|$.
Consider the $\tilde f$-$QR$ step
$F_{\tilde f}: \Aa_{\Mu,\tilde f} \to \Oo_\Mu$.
For any $S \in \Aa_{\Mu,\tilde f}$ and $r = 1, \ldots, n - 1$,
$\trace(I_r F_{\tilde f}(S)) \ge \trace(I_r S)$.
For $r = 1$, equality only holds if $(S)_{1,j} = 0$ for all $j > 1$.
\end{lemma}

This argument follows closely the first proof in \cite{DRTW}.

\smallskip

\proof
Let $V_r$ be the range of $I_r$ and
$\mu_{r,j}(S)$ be the eigenvalues
of the leading principal $r \times r$ minor of $S$,
listed in nondecreasing order.
We claim that $\mu_{r,j}(F_{\tilde f}(S)) \ge \mu_{r,j}(S)$,
which immediately implies $\trace(I_r  F_{\tilde f}(S)) \ge \trace(I_r S)$.
Recall that $F_{\tilde f}(S) = Q_\star^\ast S Q_\star$
where $Q_\star R_\star = \tilde f(S)$.
Let $U$ be an upper triangular matrix such that
$Q_\star u = \tilde f(S) U u$ for $u \in V_r$.
By min-max,
\begin{gather*} \mu_{r,j}(S) =
\max_{ \begin{matrix} \scriptstyle A \subset V_r \\ \scriptstyle \dim(A) = r+1-j \end{matrix}} \;
\min_{u \in A \smallsetminus \{0\}}
\frac{\langle u, S u \rangle}{\langle u, u \rangle}, \\
\mu_{r,j}(F_{\tilde f}(S)) =
\max_{A} \min_{u}
\frac{\langle u, F_{\tilde f}(S) u \rangle}{\langle u, u \rangle} =
\max_{A} \min_{u}
\frac{\langle \tilde f(S)U u, S \tilde f(S)U u \rangle}
{\langle \tilde f(S)U u, \tilde f(S)U u \rangle} \\
=  \max_{A' = U A} \;\min_{u' \in A' \smallsetminus \{0\} }
\frac{\langle \tilde f(S)u', S \tilde f(S)u' \rangle}
{\langle \tilde f(S)u', \tilde f(S)u' \rangle}
\end{gather*}
Notice that since $U$ is upper triangular,
the map taking $A \subset V_r$ to $A' = U A$
is a bijection among subspaces of $V_r$ of given dimension.
Since $S$ and $\tilde f(S)$ are symmetric and commute,
\[ \mu_{r,j}(F_{\tilde f}(S)) = \max_{A} \min_{u}
\frac{\langle u, S g(S) u \rangle}{\langle u, g(S) u \rangle}, \]
where $g(x) = (\tilde f(x))^2$.
The claim now follows from the inequality
\[ \langle u, u \rangle \langle u, S g(S) u \rangle -
\langle u, S u \rangle \langle u, g(S) u \rangle \ge 0. \]
Diagonalize $S= Q^\ast \Mu Q$ and
$g(S) = Q^\ast g(\Lambda)Q$ and write $Qu = (x_1,\ldots,x_n)$ so that
\[ 2 \left(\langle u, u \rangle  \langle u, S g(S) u \rangle -
\langle u, S u \rangle \langle u, g(S) u \rangle \right) =
\sum_{k,\ell} (\mu_k - \mu_\ell)
(g(\mu_k) - g(\mu_\ell)) x_k^2 x_\ell^2 \ge 0.\]
Consider now equality for the case $r = 1$.
Notice that, by hypothesis, if $k \ne \ell$ then
$(\mu_k - \mu_\ell)(g(\mu_k) - g(\mu_\ell)) > 0$.
In the max-min formula for $\trace(I_1 S) = \mu_{1,1}(S)$,
it suffices to take $u = e_1$.
Equality therefore holds only if $Qe_1$ is a canonical vector,
which implies $(S)_{1,j} = 0$ for all $j > 1$.
\qed

{\nobf Proof of Proposition \ref{prop:height}:}
For all $s \in \Ii_i$ and any
distinct eigenvalues $\lambda_j$ and $\lambda_k$,
$|\lambda_j - \lambda_i| < |\lambda_k - \lambda_i|$
if and only if $\lolo_i(\lambda_j) < \lolo_i(\lambda_k)$.
For $s \in \Ii_i$, $f(x) = x-s$, $h(x) = \lolo_i(x)$ and
$\mu_j = \lolo_i(\lambda_j)$,
define $\tilde f: \RR \to \RR$ as in Lemma \ref{lemma:ftilf}.
The function $\tilde f$ satisfies the hypothesis
of Lemma \ref{lemma:increasing}:
$\mu_j < \mu_k$ implies $|\tilde f(\mu_j)| < |\tilde f(\mu_k)|$.
Thus, by Lemma \ref{lemma:increasing},
$\trace(W  F_{\tilde f}(S)) \ge \trace(WS)$ for all $S \in \OSa$.
For $T \in \Dd^i_{\Lambda,\epsilon_{\textrm{\rm ap}}}$, take $S = h(T)$:
by Lemma \ref{lemma:ftilf}, $F_{\tilde f}(h(T)) = h(F_f(T))$
and therefore $\trace(W h(F_f(T))) \ge \trace(W h(T))$.
Again by Lemma \ref{lemma:increasing},
equality happens only if $T$ is diagonal.
Thus, $H_i$ is a height function.
Finally, choosing $\delta_H$ sufficiently small
guarantees that $H_i$ is large in $\Dd^i_{\Lambda,0}$
and small in $\partial\Dd^i_{\Lambda,\epsilon_{\textrm{\rm ap}}}$,
completing the proof.
\qed

Thus, simple shift strategies admit height functions near the deflation set.
Our reason for constructing a height function is to control the time the
sequence $(F_\sigma^k(T))$ stays in a compact set.

Assuming $\Lambda$ to be a.p.\ free,
for a shift strategy $\sigma: \ILa \to \RR$
set $\epsilon_\sigma = \epsilon_{\textrm{\rm ap}}/(1+C_\sigma)$
(where $C_\sigma$ is the constant in the definition of a simple shift strategy).
Notice that $T \in  \Dd^i_{\Lambda,\epsilon_{\sigma}}$
implies $\sigma(T) \in \Ii_i =
[\lambda_i - \epsilon_{\textrm{\rm ap}}, \lambda_i + \epsilon_{\textrm{\rm ap}}]$.

\begin{coro}
\label{coro:byebyecompact}
Let $\Lambda$ be a real diagonal $n \times n$ a.p.\ free matrix,
$\sigma$ a simple shift strategy and
$\Dd^i_{\Lambda,\epsilon_{\sigma}}$ as above.
Let $\Kk \subset \Dd^i_{\Lambda,\epsilon_{\sigma}}$ be a compact set
with no diagonal matrices:
there exists $K \in \NN$ such that
for all $T \in \Dd^i_{\Lambda,\epsilon_{\sigma}}$
there are at most $K$ points of the form $F_\sigma^k(T)$ in $\Kk$.
\end{coro}

The plan is to take $\Kk$ containing
$\calS_\sigma \cap \Dd^i_{\Lambda,\epsilon_{\sigma}}$:
the hypothesis in Theorem \ref{theo:big}
that diagonal matrices do not belong to
the singular support $\calS_\sigma$ is then natural.

\proof
Let $m_-$ be the minimum jump in $\Kk$ and
$m_+$ the size of the image of $H_i$:
\[ m_- = \inf_{T \in \Kk,\; s \in \Ii_i} H_i(F_s(T)) - H_i(T), \quad
m_+ = \sup_{T \in \Dd^i_{\Lambda,\epsilon_{\sigma}}} H_i(T) -
\inf_{T \in \Dd^i_{\Lambda,\epsilon_{\sigma}}} H_i(T). \]
By Proposition \ref{prop:height}
and the compactness of $\Kk \times \Ii_i$, $s > 0$:
take $K$ such that $Km_- > m_+$.
For a given $T$, let 
$X = \{ k \in \NN \;|\; F_\sigma^k(T) \in \Kk \}$:
we have
\[ m_+ \ge \sum_{k \in X} H_i(F_\sigma^{k+1}(T)) - H_i(F_\sigma^k(T))
\ge |X| m_- \]
and therefore $|X| < K$.
\qed

{\nobf Proof of Theorem \ref{theo:big}:}
Let $\Kk_1, \Kk_2 \subset \Dd^i_{\Lambda,\epsilon_\sigma}$
be compact sets with
$\Kk_1 \cup \Kk_2 = \Dd^i_{\Lambda,\epsilon_\sigma}$,
$\calS_\sigma \cap \Dd^i_{\Lambda,0}$ disjoint from $\Kk_1$
and with no diagonal matrices in $\Kk_2$.
By Theorem \ref{theo:squeeze},
there exists $C_{\Kk_1} > 0$ such that
$|\bbb(F_\sigma(T))| \le C_{\Kk_1} |\bbb(T)|^3$ for all $T \in \Kk_1$.
By Corollary \ref{coro:byebyecompact},
there exists $K_2 \in \NN$ such that,
given $T \in \Dd^i_{\Lambda,\epsilon_\sigma}$,
at most $K_2$ points of the form $F_\sigma^k(T)$ belong to $\Kk_2$.
In particular, there are at most $K_2$ values of $k$
for which the estimate
$|\bbb(F_\sigma^{k+1}(T))| \le C_{\Kk_1} |\bbb(F_\sigma^k(T))|^3$
does not hold.
\qed

\section{Convergence rates for a.p.\ spectra}

The aim of this section is to prove Theorem \ref{theo:bigg}.
An a.p.\ matrix $T \in \Tt$ with simple spectrum is
\textit{strong a.p.}\ if three consecutive eigenvalues
are in arithmetic progression and \textit{weak a.p.}\ otherwise.

In the a.p.\ free case discussed in the previous sections,
for an initial condition $T \in \Dd^i_{\Lambda,\epsilon}$,
the sequence $F_\sigma^k(T)$ converges to a diagonal matrix;
this follows from the fact that $\sigma(T) \approx \lambda_i$
for $T \in \Dd^i_{\Lambda,\epsilon}$.
For weak a.p.\ spectra, convergence to a diagonal matrix may not occur.

Assume $\Lambda$ to be weak a.p.
Let $\bbb_2(T) = T_{n-1,n-2}$ be the second-last subdiagonal entry;
for consistency, write $\bbb_1(T) = \bbb(T)$.
For any $i$, there exists a unique index
$c(i)$ such that $\lambda_{c(i)}$ is the eigenvalue closest to $\lambda_i$.
As we shall see, if $T \in \Dd^i_{\Lambda,\epsilon}$ then
\[ \lim_{k \to \infty} \bbb_1(F_\sigma^k(T)) = \lim_{k \to \infty} \bbb_2(F_\sigma^k(T)) = 0, \quad
\lim_{k \to \infty} (F_\sigma^k(T))_{n,n} = \lambda_i; \]
furthermore, if $T$ is unreduced then
\[ \lim_{k \to \infty} (F_\sigma^k(T))_{n-1,n-1} = \lambda_{c(i)}. \]

We begin with a technical lemma concerning the dynamics of steps $F_s$.
Item (b) is a variation of the power method argument used
to study the convergence of lower entries under $QR$ steps.

\begin{lemma}
\label{lemma:tube2}
Let $\Mu = \diag(\mu_1, \ldots, \mu_m)$ be a real diagonal
matrix with simple spectrum and $\Tt_\Mu \subset \Tt$
be the manifold of real $m \times m$ tridiagonal matrices similar to $\Mu$.
Let $I \subset \RR$ be a compact interval.
Assume that there exists $j$, $1 \le j \le m$, such that
\[ \mu_j \notin I, \quad
\max_{s \in I} |\mu_j - s| < \min_{k \ne j, s \in I} |\mu_k - s|. \]
Let $\Dd^j_{\Mu,\epsilon} \subset \Tt_\Mu$
be the $j$-th deflation neighborhood.
\begin{enumerate}[(a)]
\item{There exist $\epsilon > 0$ and $C \in (0,1)$ such that
for all $\epsilon' \in (0,\epsilon)$ and $s \in I$ we have
$F_s(\Dd^j_{\Mu,\epsilon'}) \subset \Dd^j_{\Mu,C\epsilon'}$.}
\item{Consider $T_0 \in \Tt_\Mu$ unreduced,
a sequence $(s_k)$ of elements of $I$ and $\epsilon > 0$.
Define $T_{k+1} = F_{s_k}(T_k)$.
Then there exists $k$ such that $T_k \in \Dd^j_{\Mu,\epsilon}$.}
\end{enumerate}
\end{lemma}

This will be used to study $\bbb_2(T)$ for $T \in \Dd^i_{\Lambda,\epsilon}$,
setting $I = [\lambda_i - \epsilon, \lambda_i + \epsilon]$, $j = c(i)$,
$\Mu = \Lambda_i = \diag(\lambda_1, \ldots, \lambda_{i-1},
\lambda_{i+1}, \ldots, \lambda_n)$, 
with the natural identification
between $\Tt_\Mu$ and $\Dd^i_{\Lambda,0}$.

\smallskip

\proof
Let $\tilde C \in (0,1)$ be such that
\[ \max_{s \in I} |\mu_j - s| < \tilde C \;\min_{k \ne j, s \in I} |\mu_k - s|. \]
Write
\[ r(s,T) = \frac{(R_\star)_{m,m}}{(R_\star)_{m-1,m-1}}, \quad
T - sI = Q_\star R_\star. \]
Recall from Lemma \ref{lemma:basic} and
Corollary \ref{coro:tubularQR} that $\bbb(F_s(T)) = r(s,T)\;\bbb(T)$.
We claim that for all $T \in \Dd^j_{\Mu,0}$ and $s \in I$, $|r(s,T)| \le \tilde C$.
Since $T \in \Dd^j_{\Mu,0}$, $|(R_\star)_{m,m}| = |\mu_j - s|$.
Let $R_-$ be the leading principal minor of $R_\star$ of order $m-1$:
its singular values are $|\mu_k - s|$, $k \ne s$.
In particular, all singular values
are larger that $|(R_\star)_{m,m}|/{\tilde C}$.
Thus
\[ |(R_\star)_{m-1,m-1}| = \| e_{m-1}^\ast R_- \| \ge
\frac{|(R_\star)_{m,m}|}{\tilde C} \| e_{m-1} \| =
\frac{|(R_\star)_{m,m}|}{\tilde C}, \]
proving our claim.
Take  $C = (1+\tilde C)/2$:
by continuity, for sufficiently small $\epsilon > 0$, we have $|r(s,T)| < C$
for all $T \in \Dd^j_{\Mu,\epsilon}$, $s \in I$.
Thus, for $T \in \Dd^j_{\Mu,\epsilon}$ and $s \in I$,
$|\bbb(F_s(T))| \le C\,|\bbb(T)|$;
item (a) follows.


For item (b), write $T_{k+1} = Q_k^\ast T_k Q_k$
where $T_k - s_k I = Q_k R_k$ is a $Q_\star R_\star$ decomposition.
Notice that, by hypothesis, $I$ is disjoint from the spectrum so that
$T_0 - s_0 I$ is invertible.
We have $(T_0 - s_0 I)^{-1} = R^{-1} Q_0^\star$
so the rows of $Q_0^\star$ are obtained from those of $(T_0 - s_0 I)^{-1}$
by Gram-Schmidt from bottom to top.
In particular, $Q_0 e_m = c_0 (T_0 - s_0 I)^{-1} e_m$,
$c_0 > 0$.
More generally, we claim that
\begin{gather*}
P_k e_m =
c (T_0 - s_{k-1} I)^{-1} \cdots (T_0 - s_1 I)^{-1} (T_0 - s_0 I)^{-1} e_m, \\
c > 0, \quad P_k = Q_0 Q_1 \cdots Q_{k-1} \in SO(m).
\end{gather*}
Indeed, by induction and using that $T_1 = Q_0^\ast T_0 Q_0$,
\begin{align*}
P_k e_m &= c' Q_0 (T_1 - s_{k-1} I)^{-1} \cdots (T_1 - s_1 I)^{-1} e_m  \\
&= c' (T_0 - s_{k-1} I)^{-1} \cdots (T_0 - s_1 I)^{-1} Q_0 e_m \\
&= c (T_0 - s_{k-1} I)^{-1} \cdots (T_0 - s_1 I)^{-1} (T_0 - s_0 I)^{-1} e_m.
\end{align*}
(Integrability of the Toda lattice is present here yet another time.)
For $\alpha = 1, \ldots, m$,
let $v_\alpha$ be the unit eigenvector associated to $\mu_\alpha$.
We claim that
\[ \lim_{k \to \infty} P_k e_m = \pm v_j. \]
Indeed, write $e_m = \sum_{\alpha=1}^m a_\alpha v_\alpha$,
where $a_\alpha = \langle v_\alpha, e_m \rangle$ is the last
coordinate of $v_\alpha$.
It is well known that the last coordinates of the eigenvectors $v_\alpha$
of the unreduced matrix $T$ are nonzero:
in particular, $a_j \ne 0$; assume without loss $a_j > 0$.
We have
\begin{align*}
P_k e_m &=
c (T_0 - s_{k-1} I)^{-1} \cdots (T_0 - s_1 I)^{-1} (T_0 - s_0 I)^{-1} e_m \\
&= c \sum_{\alpha = 1}^m
\frac{a_\alpha}{(\mu_\alpha - s_{k-1})\cdots(\mu_\alpha - s_{0})}
v_\alpha
= c_k \left( v_j + \sum_{\alpha \ne j} b_{k,\alpha} v_\alpha \right),\\
&\qquad c_k > 0,  \quad
b_{k,\alpha} =
\frac{a_\alpha}{a_j}\;\frac{\mu_j - s_{k-1}}{\mu_\alpha - s_{k-1}}\cdots
\frac{\mu_j - s_{0}}{\mu_\alpha - s_{0}}.
\end{align*}
Since $| \mu_j - s_{k-1} |/|\mu_\alpha - s_{k-1}| < \tilde C$
we have $|b_{k,\alpha}| \le (\tilde C)^k\,|a_\alpha/a_j|$ and
therefore $\lim_{k \to \infty} b_{k,\alpha} = 0$,
proving the claim. We have
\begin{gather*}
\lim_{k \to \infty} \bbb(T_k) = \lim_{k \to \infty} (T_k)_{m,m-1} =
\lim_{k \to \infty} e_{n-1}^\ast T_k e_m =
\lim_{k \to \infty} (P_k e_{m-1})^\ast T_0 (P_k e_m)= \\
= \lim_{k \to \infty} (P_k e_{m-1})^\ast \mu_j (P_k e_m) +
\lim_{k \to \infty} (P_k e_{m-1})^\ast (T_0 - \mu_j I) (P_k e_m).
\end{gather*}
The first limit in the last expression is zero because
$P_k e_{m-1} \perp P_k e_{m}$;
the second is zero because $P_k e_{m-1}$ is bounded
and
\[ \lim_{k \to \infty} (T_0 - \mu_j I) (P_k e_m)
= (T_0 - \mu_j I) \lim_{k \to \infty}  (P_k e_m)
= (T_0 - \mu_j I) v_j = 0. \]
\qed

Consider the \textit{double deflation set}
$\Cc_{\Lambda,0} \subset \Dd_{\Lambda,0} \subset \Tt_\Lambda$:
\[ \Cc_{\Lambda,0} = \{ T \in \Tt_{\Lambda} \;|\; \bbb_1(T) =
\bbb_2(T) = 0 \}. \]
For Wilkinson's strategy $\omega$, it turns out that the set $\Cc_{\Lambda,0}$
is disjoint from the singular support $\calS_\omega$.  More generally, if a
shift strategy $\sigma$ satisfies $\Cc_{\Lambda,0} \cap \calS_\sigma =
\emptyset$ then cubic convergence of $F_\sigma$ holds even for weak a.p.\
spectra: this is Theorem \ref{theo:bigg}, which we prove below.

In \cite{LST2}, we show examples of
unreduced tridiagonal $3 \times 3$ matrices
with spectrum $-1, 0, 1$ for which Wilkinson's shift $F_\omega$ converges
quadratically to a reduced but not diagonal matrix
in the singular support $\calS_\omega$.
Similarly, we conjecture that
for strong a.p.\ diagonal $n \times n$ matrices $\Lambda$
there exists a set $\Xx \subset \ILa$ of Hausdorff codimension 1
of unreduced matrices $T$ for which $F_\omega^k(T)$ converges quadratically
to a matrix in $\calS_\omega \cap \Dd_{\Lambda,0}$ with $T_{n-1,n-2} \ne 0$.

%
%

With the natural identification between $\Dd^i_{\Lambda,0}$ and $\Tt_{\Lambda_i}$,
we may consider $\Dd^j_{\Lambda_i,\epsilon_2}$ to be a subset of $\Dd^i_{\Lambda,0}$.
Let
\[ \Cc^{j,i}_{\Lambda,\epsilon_2,\epsilon_1} =
\Dd^i_{\Lambda,\epsilon_1} \cap
\Pi_i^{-1}(\Dd^{j}_{\Lambda_i,\epsilon_2}). \]
For small $\epsilon_1, \epsilon_2 > 0$,
$T \in \Cc^{j,i}_{\Lambda,\epsilon_2,\epsilon_1}$
implies
\[ T_{n-1,n-1} \approx \lambda_j, \quad T_{n,n} \approx \lambda_i,
\quad \bbb_1(T) \le \epsilon_1, \quad \bbb_2(T) \approx 0. \]
These compact sets turn out to be manifolds with corners
but we shall neither prove nor use this fact.
Lemma \ref{lemma:tube2} can be rephrased in terms of the sets $\Cc^{j,i}_{\Lambda,\epsilon_2,\epsilon_1}$.

\begin{coro}
\label{coro:tube3}
Let $\Lambda$ to be weak a.p. spectrum
and $\sigma$ be a simple shift strategy.
There exists $\epsilon > 0$ such that,
for all $i$ and for all $\epsilon_1 \in (0,\epsilon)$:
\begin{enumerate}[(a)]
\item{there exists $C \in (0,1)$ such that,
for all sufficiently small $\epsilon_2 > 0$
we have $F_\sigma(\Cc^{c(i),i}_{\Lambda,\epsilon_2,\epsilon_1}) \subset
\Cc^{c(i),i}_{\Lambda,C \epsilon_2,\epsilon_1}$;}
\item{for all unreduced $T \in \Dd^i_{\Lambda,\epsilon}$
and for all $\epsilon_1, \epsilon_2 > 0$ there exists $k$ such that
$F_\sigma^k(T) \in \Cc^{c(i),i}_{\Lambda,\epsilon_2,\epsilon_1}$.}
\end{enumerate}
\end{coro}

\proof
Combine Lemma \ref{lemma:tube2} with
$\Pi_i \circ F_s = F_s \circ \Pi_i$
(Proposition \ref{prop:pia}).
\qed

{\nobf Proof of Theorem \ref{theo:bigg}:}
From the hypothesis that $\Cc_{\Lambda,0}$ and $\calS_\sigma$ are disjoint
it follows that, for sufficiently small $\epsilon_1, \epsilon_2 > 0$,
the shift strategy $\sigma$ is smooth in $\Cc^{c(i),i}_{\Lambda,\epsilon_2,\epsilon_1}$.
As in Lemma \ref{lemma:wlip}, from a Taylor expansion around $T_0 \in \Dd^i_{\Lambda,0}$,
there exists $C_2$ such that $|\sigma(T)| \le C_2 |\bbb_1(T)|^2$ for all
$T \in \Cc^{c(i),i}_{\Lambda,\epsilon_2,\epsilon_1}$.
As in the proof of Theorem \ref{theo:squeeze}, there exists $C_3$
such that $|\bbb_1(F_\sigma(T))| \le C_3 |\bbb_1(T)|^3$
for all $T \in \Cc^{c(i),i}_{\Lambda,\epsilon_2,\epsilon_1}$.
From item (a) of Corollary \ref{coro:tube3},
$\Cc^{c(i),i}_{\Lambda,\epsilon_2,\epsilon_1}$
is invariant under $F_\sigma$;
from item (b), for all unreduced $T \in \Dd^i_{\Lambda,\epsilon}$
(where $\epsilon$ is sufficiently small)
there exists $K$ such that, for all $k > K$, $F_\sigma^k(T) \in \Cc^{c(i),i}_{\Lambda,\epsilon_2,\epsilon_1}$,
completing the proof.
\qed

\section{Two counterexamples}

In this section we present two examples which show that
natural strengthenings of Theorems \ref{theo:big} and \ref{theo:bigg}
do not hold for Wilkinson's strategy $\omega$.

\begin{figure}[ht]
\begin{center}
\psfrag{Yy}{$\calS_\omega$}
\psfrag{S}{$S$}
\epsfig{height=28mm,file=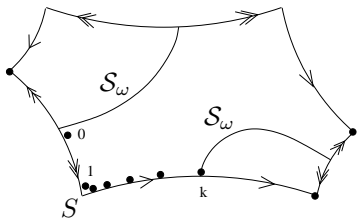}
\end{center}
\caption{\capsize We may have $F_\omega^k(T) \in \calS_\omega$ for large values of $k$.}
\label{fig:notsocubic}
\end{figure}

We use the notation of Section 3.
In Figure \ref{fig:notsocubic}, where $\Lambda = \diag(1,2,4)$,
we indicate a sequence $\tilde F_\omega^k(T)$ 
which enters the deflation neighborhood $\Dd^i_{\Lambda,\epsilon}$
near one diagonal matrix but travels within the neighborhood
towards another diagonal matrix.
Theorem \ref{theo:squeeze} guarantees the cubic decay of the $(3,2)$
entry whenever $\tilde F_\omega^k(T)$ stays away from the singular support $\calS_\omega$.
Consistently with Theorem \ref{theo:big},
this happens for practically all values of $k$.
Notice however that no uniform bound exists on the number of iterations
needed to reach (a neighborhood of) $\calS_\omega$.
As proved in \cite{LST2}, in this instance cubic decay does not hold.
More precisely,
it is {\it not} true that given an a.p.\ free matrix $\Lambda$
there exist $C > 0$ and $K$ such that
$|\bbb(F_\omega^{k+1}(T))| \le C |\bbb(F_\omega^k(T))|^3$ for all $k > K$.

Consider now the weak a.p. spectrum $\Lambda = \diag(-1,0,0.3,1)$ and
\[ T_0 = \begin{pmatrix} 0.3 & 0 \\ 0 & S_0 \end{pmatrix} \in \ILa \]
where $S_0 \in \Tt_{\Lambda_3}$, $\Lambda_3 = \diag(-1,0,1)$,
is an example of unreduced matrix obtained in \cite{LST2}
for which convergence is strictly quadratic, i.e.,
\[ C_- |\bbb(F_\omega^k(S_0))|^2 < |\bbb(F_\omega^{k+1}(S_0))| <
C_+ |\bbb(F_\omega^k(S_0))|^2, \]
for all $k$, where $0 < C_- < C_+$.
Trivially, the analogous estimate holds for $\bbb(F_\omega^k(T_0))$.
By sheer continuity, given $K$, there exists $\epsilon > 0$ such that if
$T \in \ILa$ satisfies $\| T - T_0 \| < \epsilon$  then
\[ C_- |\bbb(F_\omega^k(T))|^2 < |\bbb(F_\omega^{k+1}(T))| <
C_+ |\bbb(F_\omega^k(T))|^2 \]
still holds for all $k < K$.
Thus, the uniform estimate in Theorem \ref{theo:big}
fails for weak a.p. spectra, even for unreduced matrices.

\section{Appendix: Proof of Theorem \ref{theo:connect} }

Recall from Section 2 that $E_j \in \cE$ is defined by
\[ (E_j)_{i,i} = \begin{cases} 1, & i \le j, \\ -1, & i > j; \end{cases} \]
the involutions $\eta_j$ are defined by $\eta_j(T) = E_j T E_j$,
which differs from $T$ only in the sign of the $j$-th subdiagonal coordinate.
Let $\cM_j \subset \ILa$ be the \textit{mirror}, i.e.,
the set of fixed points of $\eta_j$:
for $T \in \ILa$ we have $T \in \cM_j$ if and only if $(T)_{j+1,j} = 0$.
Let $S_n$ be the symmetric group of permutations $\pi$
of the set $\{1, 2, \ldots, n\}$.
For $\pi \in S_n$, let $\cM_{j,\pi} \subset \cM_j$
be the set of matrices for which 
the eigenvalues of the top $j \times j$ principal
subblock are $\lambda_{\pi(1)}, \ldots, \lambda_{\pi(j)}$ so that
\[ \cM_{j,\pi} \approx
\cT_{\diag(\lambda_{\pi(1)}, \ldots, \lambda_{\pi(j)})} \times
\cT_{\diag(\lambda_{\pi(j+1)}, \ldots, \lambda_{\pi(n)})}. \]
Thus, $\cM_j$ is a submanifold of codimension $1$ with
$\binom{n}{j}$ connected components $\cM_{j,\pi}$.

The diagonal matrices in $\ILa$ are labeled by $\pi \in S_n$: let 
\[ \Lambda^\pi = \diag(\lambda_{\pi(1)}, \ldots, \lambda_{\pi(n)}). \]
Let $\JLa \subset \ILa$ be the set of tridiagonal matrices
with nonnegative subdiagonal entries.
The set $\JLa$ is homeomorphic to the \textit{permutohedron}
$\cP_\Lambda$ (\cite{Tomei}),
the convex hull of the points 
\[ v_\pi = (\lambda_{\pi^{-1}(1)}, \ldots, \lambda_{\pi^{-1}(n)}) \in \RR^n,
\qquad \pi \in S_n; \]
the vertices of $\cP_\Lambda$ are $v_\pi$.
An explicit homeomorphism takes $T = Q^\ast \Lambda Q$ to
the vector in $\RR^n$ whose $j$-th coordinate is
$(Q \Lambda Q^\ast)_{jj}$;
this map takes $\Lambda^\pi$ to $v_\pi$ (\cite{BFR}).
We use this map to endow $\JLa$ with a combinatorial structure
of vertices, faces and hyperfaces:
in particular, vertices of $\JLa$ are diagonal matrices.
It turns out that the hyperfaces of $\JLa$
are the intersections $\cM_{j,\pi} \cap \JLa$.

\bigskip

\begin{lemma}
\label{lemma:polytope}
Let $\cP \subset \RR^n$ be the convex hull of a finite set.
Let $\tilde F: \cP \to \cP$ be a continuous function.
Assume that for any hyperface $\cQ \subset \cP$
we have $\tilde F(\cQ) \subset \cQ$.
Then $\tilde F$ is surjective.
\end{lemma}

\proof
The dimension of a convex subset of $\RR^n$
is the dimension of the affine subspace
spanned by its vertices.
Notice that any face (of any dimension)
is the intersection of hyperfaces and therefore also
invariant under $\tilde F$.

We use relative homology: if the dimension of $\cP$ is $d$
then $H_d(\cP,\partial \cP) = \ZZ$;
we prove that $\tilde F_\ast: H_d(\cP,\partial \cP) \to H_d(\cP,\partial \cP)$
is the identity.
This implies the lemma: if $x_0$ is an interior point of $\cP$
not in the image of $\tilde F$ 
then since $H_d(\cP,\cP \smallsetminus \{x_0\}) = H_d(\cP,\partial \cP)$
we have $\tilde F_\ast = 0$, a contradiction.

The proof of the claim is by induction on the dimension $d$ of $\cP$.
The case $d = 0$ is trivial;
in the case $d = 1$ the polytope $\cP$ is an interval
and $\tilde F$ takes each endpoint to itself
and again the claim is easy.
In general, let $Q$ be a hyperface of $\cP$
so that the dimension of $\cQ$ is $d-1$ and, by induction,
$\tilde F_\ast: H_{d-1}(\cQ,\partial \cQ) \to H_{d-1}(\cQ,\partial \cQ)$
is the identity.
We have $\cQ \subset \partial \cP$
and $H_{d-1}(\cQ,\partial \cQ) =
H_{d-1}(\partial \cP, \partial \cQ \cup (\partial \cP \smallsetminus \cQ)) =
H_{d-1}(\partial \cP)$ and therefore
$\tilde F_\ast: H_{d-1}(\partial \cP) \to H_{d-1}(\partial \cP)$
is the identity.
Since $\cP$ is contractible, the long exact sequence for
relative homology implies that
$\tilde F_\ast: H_d(\cP,\partial \cP) \to H_d(\cP,\partial \cP)$
is the identity, completing the proof.
\qed

{\noindent\bf Proof of Theorem \ref{theo:connect}:}
For (a), first notice that the condition
$F \circ \eta_i = \eta_i \circ F$ implies $F(\cM_i) \subseteq \cM_i$.
Since diagonal matrices are fixed points this implies
$F(\cM_{i,g}) \subseteq \cM_{i,g}$.
Restrict $F$ to $\JLa$ and drop signs to define a continuous
map $\tilde F: \JLa \to \JLa$ which keeps
each hyperface of $\JLa$ invariant.
By Lemma \ref{lemma:polytope},
$\tilde F$ is surjective and therefore (by equivariance) so is $F$.

For (b), let $\Bb^i \subset \ILa$ be
the basins of attraction of each invariant neighborhood
$\interior(\cK_i)$, i.e., 
$T \in \Bb^i$ if there exists $k \in \NN$
such that $F_\sigma^k(T) \in \interior(\cK_i)$.
The sets $\Bb^i$ are clearly disjoint with
$\cK_i \subset \Bb^i$.
They are also open subsets of $\ILa$
since $\Bb_i = \bigcup_k F^{-k}(\interior(\cK_i))$.
Since $\ILa$ is connected there exists $T \notin \bigcup_i \Bb^i$
and we are done.
\qed

\bigskip\bigskip\bigbreak

{

\parindent=0pt
\parskip=0pt
\obeylines

Ricardo S. Leite, Departamento de Matemática, UFES
Av. Fernando Ferrari, 514, Vitória, ES 29075-910, Brazil

\smallskip

Nicolau C. Saldanha and Carlos Tomei, Departamento de Matem\'atica, PUC-Rio
R. Marqu\^es de S. Vicente 225, Rio de Janeiro, RJ 22453-900, Brazil

\smallskip
rsleite@pq.cnpq.br
saldanha@puc-rio.br; http://www.mat.puc-rio.br/$\sim$nicolau/
tomei@mat.puc-rio.br

}

\end{document}